\documentclass[12pt,thmsa]{article}
\usepackage{amssymb}
\usepackage{amsfonts}
\usepackage{amsmath}

\begin{document}

\begin{center}
\bigskip

\bigskip\ {\LARGE Modules with ascending chain condition on annihilators and
Goldie modules}

\bigskip

{\Large Jaime Castro P\'{e}rez}

Escuela de Ingenier\'{\i}a y Ciencias

Instituto Tecnol\'{o}gico y de Estudios Superiores de Monterrey

Calle del Puente 222, Tlalpan

14380 M\'{e}xico, D.F.

M\'{e}xico

\bigskip

{\Large \ Mauricio Medina B\'{a}rcenas and Jos\'{e} R\'{\i}os Montes}

Instituto de Matem\'{a}ticas

Universidad Nacional Aut\'{o}noma de M\'{e}xico

Area de la Investigaci\'{o}n Cient\'{\i}fica

Circuito Exterior, C.U.

04510 M\'{e}xico, D.F.

M\'{e}xico

\bigskip

\textbf{Abstract}
\end{center}

\begin{quotation}
Using the concepts of prime module, semiprime module and the concept of
ascending chain condition (ACC) on annihilators for an $R$-module $M$ . We
prove that if \ $M$ is semiprime \ and projective in $\sigma \left[ M\right]
$, such that $M$ satisfies ACC on annihilators, then $M$ has finitely many
minimal prime submodules. Moreover if each submodule $N\subseteq M$ contains
a uniform submodule, we prove that there is a bijective correspondence
between a complete set of representatives of isomorphism classes of
indecomposable non $M$-singular injective modules in $\sigma \left[ M\right]
$ and the set of minimal primes in $M$. If $M$ is Goldie module then $%
\widehat{M}\cong E_{1}^{k_{1}}\oplus E_{2}^{k_{2}}\oplus ...\oplus
E_{n}^{k_{n}}$ where each $E_{i}$ is a uniform $M$-injective module. As an
application, new characterizations of left Goldie rings are obtained.

\bigskip


\bigskip


\bigskip
\end{quotation}

\section*{Introduction}

\bigskip

The Rings with ascending chain condition (ACC) on annihilators has been
studied by several authors. For non commutative semiprime rings A.W.
Chatters and C. R. Hjarnavis in [6] showed that if $R$ is semiprime ring
such that $R$ satisfies ACC on annihilator ideals, then $R$ has finitely
many minimal prime ideals. K.R. Goodearl and R. B. Warfield in [8] showed
that if $R$ is a semiprime right Goldie ring, then $R$ has finitely many
minimal prime ideals.

In this paper we work in a more general situation, using the concept of
prime and semiprime submodules defined in [11], [13] and the concept of
ascending chain condition (ACC) on annihilators for a $R$-module $M$ defined
in [4], we show that if $M$ is projective in $\sigma \left[ M\right] $ and
semiprime module, such that $M$ satisfies ACC on annihilators, then $M$ has
finitely many minimal prime submodules in $M$. Moreover if each submodule $%
N\subseteq M$ contains a uniform submodule, we prove that there is a
bijective correspondence between a complete set of isomorphism classes of
indecomposable non $M$-singular injective modules in $\sigma \left[ M\right]
$ $\mathcal{E}_{\chi \left( M\right) }(M)$ and the set of minimal prime in $%
M $ submodules $\ Spec^{Min}\left( M\right) $. Also we prove that there is a
bijective correspondence between sets $\mathcal{E}_{\chi \left( M\right)
}(M) $ and $P_{M}^{Min}$, where $P_{M}^{Min}=\left\{ \chi \left( M/P\right)
\mid P\subseteq M\text{ is a minimal prime in }M\right\} $. So we obtain new
information about semiprime rings with ACC on annihilators. As an
application, of these results we obtain a new characterization of Goldie
Rings in terms of the of minimal prime idelas of $R$ and a new
characterization of these Rings in terms of the continuous modules with ACC
on annihilators.

In order to do this, we organized the article in three sections. Section 1,
we give some results about prime and semiprime modules and we define the
concept right annihilator for a submodule $N$ of $M$. In particular we prove
in Proposition 1.16, that if $M$ is projective in $\sigma \left[ M\right] $
and $M$ is semiprime module, then $Ann_{M}\left( N\right) =Ann_{M}^{r}\left(
N\right) $ for all $N\subseteq M$. In Section 2, we define the concept of
ascending chain condition (ACC) on annihilators for an $R$-module $M$ and we
give the main result of this paper Theorem 2.2, we prove that if $M$ is a
semiprie module projective in $\sigma \left[ M\right] $ and $M$ satisfies
ACC on annihilators, then $M$ has finitely many minimal prime submodules in $%
M$. Moreover if $P_{1}$, $P_{2},...,P_{n}$ are the minimal prime in $M$
submodules, then $P_{1}\cap P_{2}\cap $... $\cap P_{n}=0$. In Theorem 2.7,
we show that there is a bijective correspondence between sets $\mathcal{E}%
_{\chi \left( M\right) }(M)$ and $Spec^{Min}\left( M\right) $. Moreover in
Corollary 2.13 we obtain an application for semiprime rings. We prove that
if $R$ is a semiprime ring such $R$ satisfies ACC on left annihilator and
each $0\neq I\subseteq R$ left ideal of $R$, contains a uniform left ideal
of $R$, then there is a bijective correspondence between the set of
representatives of isomorphism classes of indecomposable non singular
injective $R$-modules and the set of minimal prime ideals of $R$. In
particular in Corollary 2.15 we show, that for $R$ a semiprime left Goldie
ring the bijective correspondence above mentioned is true. Moreover in this
Corollary we show that $R$ has local Gabriel correspondence with respect to $%
\chi \left( R\right) =\tau _{g}$. Other important results are the Theorem
2.18 and 2.20. In the first Theorem we show that, if $M$ satisfies ACC on
annihilators and for each $0\neq N\subseteq M$, $N$ contains a uniform
submodule, then $\widehat{M}$ $=\widehat{N_{1}}\oplus \widehat{N_{2}}\oplus
...\oplus \widehat{N_{n}}$ where $N_{i}=Ann_{M}^{r}\left( P_{i}\right) $ for
$1\leq i\leq n$. And in the second Theorem we prove that if $M$

is semiprime Goldie module, then there are $E_{1}$, $E_{2}$, ..., $E_{n}$
uniform injective modules such that $\widehat{M}\cong E_{1}^{k_{1}}\oplus
E_{2}^{k_{2}}\oplus ...\oplus E_{n}^{k_{n}}$ and $Ass_{M}\left( E_{i}\right)
=\left\{ P_{i}\right\} $. As other application of these results we give the
Theorem 2.27, where we prove that if $R$ is a semiprime ring, such that $R$
satisfies ACC on left annihilators and for each non zero left ideal $%
I\subseteq R$, $I$ contains a uniform left ideal, then the following
conditions are equivalents: $i)$ $R$ is left Goldie ring, $ii)$ $P_{i}$ has
finite uniform dimension for each $1\leq i\leq n$. In the Section 3 we use
the concept of continuous module and we show in the Theorem 3.6 that if $M$
is continuous, retractable, non $M$-singular module and \ $M$ satisfies ACC
on annihilators, then $M$ is a semiprime Goldie module. Finally as one more
application we obtain the Corollary 3.7, where we show that, if $R$ is a
continuous, non singular ring and $R$ satisfies ACC on left annihilators,
then $R$ is semiprime left Goldie ring.

In what follows, $R$ will denote an associative ring with unity and $R$-$Mod$
will denote the category of unitary left $R$-modules. Let $M$ and $X$ be $R$%
-modules. $X$ is said to be $M$-generated if there exists an $R$-epimorphism
from a direct sum of copies of $M$ onto $X$. The category $\sigma \left[ M%
\right] $ is defined as the full subcategory of $R$-$Mod$ containing all $R$%
-modules $X$ which are isomorphic to a submodule of an $M$-generated module.

Let $M$-$tors$ be the frame of all hereditary torsion theories on $\sigma %
\left[ M\right] $. For a family $\{M_{\alpha }\}$ of left $R$-modules in $%
\sigma \left[ M\right] $, let $\chi \left( \{M_{\alpha }\}\right) $ be the
greatest element of $M$-$tors$ for which all the $M_{\alpha }$ are torsion
free, and let $\xi \left( \{M_{\alpha }\}\right) $ denote the least element
of $M$-$tors$ for which all the $M_{\alpha }$ are torsion. $\chi \left(
\{M_{\alpha }\}\right) $ is called the torsion theory cogenerated by the
family $\{M_{\alpha }\},$ and $\xi \left( \{M_{\alpha }\}\right) $ is the
torsion theory generated by the family $\{M_{\alpha }\}.$ In particular, the
greatest element of $M$-$tors$ is denoted by $\chi $ and the least element
of $M$-$tors$ is denoted by $\xi .$ If $\tau $ is an element of $M$-$tors$, $%
gen\left( \tau \right) $ denotes the interval $\left[ \tau ,\chi \right] $.

Let $\tau \in $ $M$-$tors$. By $\mathbb{T}_{\tau },\mathbb{F}_{\tau
},t_{\tau }$ we denote the torsion class, the torsion free class and the
torsion functor associated to $\tau $, respectively. For $N\in $ $\sigma %
\left[ M\right] $, $N$ is called $\tau $-cocritical if $N\in $ $\mathbb{F}%
_{\tau }$ and for all $0\neq N^{\prime }\subseteq N,$ $N/N^{\prime }\in
\mathbb{T}_{\tau }$. We say that $N$ is cocritical if $N$ is $\tau $%
-cocritical for some $\tau \in $ $M$-$tors$. If $N$ is an essential
submodule of $M$, we write $N\subseteq _{ess}M$. If $N$ is a fully invariant
submodule of $M$ we write $N\subseteq _{FI}M$. For $\tau \in M$-tors and $%
M^{\prime }\in \sigma \left[ M\right] $, a submodule $N$ of $M^{\prime }$ is
$\tau $-pure in $M^{\prime }$ if $M^{\prime }/N\in \mathbb{F}_{\tau }$.

A module $N\in \sigma \left[ M\right] $ is called singular in $\sigma \left[
M\right] $, or $M$-singular, if there is an exact sequence in $\sigma \left[
M\right] $ \ $0\rightarrow K\rightarrow L\rightarrow N\rightarrow 0$, with $%
K $ essential in $L$. The class $\mathcal{S}$ of all $M$-singular modules in
$\sigma \left[ M\right] $ is closed by taking submodules, factor modules and
direct sums. Therefore, any $L\in \sigma \left[ M\right] $ has a largest $M$%
-singular submodule $\mathcal{Z}\left( L\right) =\sum \left\{ f\left(
N\right) \mid N\in \mathcal{S}\text{ and }f\in Hom_{R}\left( N,L\right)
\right\} $. \ $L$ is called non $M$-singular if $\mathcal{Z}\left( L\right)
=0$. If $N$ is a fully invarint submodule of $M$, we write $N\subseteq
_{FI}M $.

Injective modules and injective hulls exist in $\sigma \left[ M\right] $,
since injective hull of $X$ in $\sigma \left[ M\right] $ is $\widehat{X}%
=tr^{M}\left( E\left( X\right) \right) =\sum_{f\in Hom_{R}\left( M,E\left(
X\right) \right) }f\left( M\right) $, where $E\left( X\right) $ denotes the
injective hull of $X$ in $R$-$Mod$.

Let $M\in R$-$Mod$. In [\textbf{1,} Definition\thinspace 1.1] the
annihilator in $M$ of a class $\mathcal{C}$ of \ $R$-modules is defined as $%
Ann_{M}(\mathcal{C})=\underset{K\in \Omega }{\cap }K$, where

$\Omega =\left\{ K\subseteq M\mid \text{there exists }W\in \mathcal{C}\text{
and }f\in Hom_{R}(M,W)\text{ with }K=\ker f\right\} $. Also in [1, \
Definition\thinspace 1.5] a product is defined in the following way. Let $%
N\,\;$\ be a submodule of $M$. For each module $X\in R$-$Mod$, $N\cdot
X=Ann_{X}(\mathcal{C})$, where $\mathcal{C}$ \ is the class of modules $W$,
such that $f(N)=0$ for all $f\in Hom_{R}(M,W)$. For $M\in R$-$Mod$ and $K$, $%
L$ submodules of $M$, in [2] is defined the product $K_{M}L$ as $%
K_{M}L=\sum \left\{ f\left( K\right) \mid f\in Hom_{R}\left( M,L\right)
\right\} $. We in [3,Proposition 1.9] that if $M\in R$-Mod and $\mathcal{C}$
be a class of left $R$-modules, then $Ann_{M}(\mathcal{C})=\sum \left\{
N\subset M\mid N_{M}X=0\,\text{\ for all }\,X\in \mathcal{C}\right\} $.

Moreover Beachy showed in [1, Proposition 5.5] that, if $M$ is projective in
$\sigma \left[ M\right] $ and $N$ is a any submodule of $M$, then $N\cdot
X=N_{M}X$, for any $R$-module $X\in \sigma \left[ M\right] $.

Let $M\in R$-$Mod$ and $N\neq M$ a fully invariant submodule of $M$, $N$ is
prime in $M$ if for any $K$, $L$ fully invariant submodules of $M$ we have
that $K_{M}L\subseteq N$ implies that $K\subseteq N$ or $L\subseteq N$. We
say that $M$ is a prime module if $0$ is prime in $M$ see [11, Definition 13
and Definition 16].

For $N\in \sigma \lbrack M]$, a proper fully invariant submodule $K$ of $M$
is said to be associated to $N$, if there exists a non-zero submodule $L$ of
$N$ such that $Ann_{M}(L^{\prime })=K$ for all non-zero submodules $%
L^{\prime }$ of $L$. By [3,\ Proposition 1.16] we have that $K$ is prime in $%
M$. We denote by $Ass_{M}\left( N\right) $ the set of all submodules prime
in $M$ associated to $N$. Also note that, if $N$ is a uniform module, then $%
Ass_{M}(N)$ has at most one element.

A module $M$ is retractable if $Hom_{R}\left( M,N\right) \neq 0$ for all $%
0\neq N\subseteq M$.

For details about concepts and terminology concerning torsion theories in $%
\sigma \left[ M\right] $, see [16 ] and [17].

\bigskip

\section{Preliminaries}

\begin{center}
\bigskip
\end{center}

The following definition was given in [4, Definition 3.2] we include it for
the convenience of the reader.

\bigskip

\bigskip \textbf{Definition 1.1.} Let $M\in R$-Mod and $N\neq M$ a fully
invariant submodule of $M$. We say that $N$ is semiprime in $M$ if for any $%
K\subseteq _{FI}M$ such that $K_{M}K$ $\subseteq N$, then $K\subseteq N$. We
say $M$ is semiprime if $0$ is semiprime in $M$.

\bigskip

\textbf{Remark 1.2}. Notice that if $M$ is projective in $\ \sigma \left[ M%
\right] $ and $N\subseteq _{FI}M$, then by [4, Remark 3.3], $N$ is semiprime
in $M$ if and only if for any submodule $K$ of $M$ such that $%
K_{M}K\subseteq N$ we have that $K\subseteq N$. Analogously by [4 Remark
3.22] we have that if $P\subseteq _{FI}M$, then $P$ is prime in $M$ if and
only if for any $K$, $L$ submodules of $M$ such that $K_{M}L\subseteq P$ we
have that $K\subseteq P$ or $L\subseteq P$.

\bigskip

\textbf{Proposition 1.3}. Let $M\in R$-Mod be projective in $\sigma \left[ M%
\right] $. If $N\subseteq _{FI}M$, then the following conditions are
equivalent

$i)$ $N$ is semiprime in $M$

$ii)$ For any submodule $K$ of $M$ such that $N\subseteq K$ and $%
K_{M}K\subseteq N$, then $K=N$.

\bigskip

\textbf{Proof}. $i)\Rightarrow ii)$ It is clear by Remark 1.2.

\bigskip

$ii)\Rightarrow i)$ $\ $Let $K\subseteq M$ a submodule such that $%
K_{M}K\subseteq N$. We claim that $\left( K+N\right) _{M}K\subseteq N$ and $%
K_{M}\left( K+N\right) \subseteq N$. In fact by [3, Proposition 1.3] we have
that $\left( K+N\right) _{M}K=K_{M}K+N_{M}K$. As $N$ is a fully invariant
submodule of $M$, then $N_{M}K\subseteq N$. Thus $\left( K+N\right)
_{M}K\subseteq N$. Now as $K_{M}K\subseteq N$, then $K_{M}K\subseteq N\cap K$%
. So by [1, Proposition 5.5] we have that $K_{M}\left( \dfrac{K}{N\cap K}%
\right) =0$. Hence $K_{M}\left( \dfrac{K+N}{N}\right) =0$. Newly by [\textbf{%
1,} Proposition 5.5] $K_{M}\left( K+N\right) \subseteq N$. Hence by [3,
Proposition 1.3] we have that $\left( K+N\right) _{M}\left( K+N\right)
=K_{M}\left( K+N\right) +N_{M}\left( K+N\right) \subseteq N$. By hypothesis
we have that $\left( K+N\right) =N$. Thus $K=N$.

\bigskip

\textbf{Remark 1.4}. Let \ $M\in R$-Mod be projective in $\sigma \left[ M%
\right] $. Similarly to the proof of Proposition 1.3, we can prove that a
fully invariant submodule $P$ of $M$ is \ prime in $M$ if and only if for
any submodules $K$ and $L$ of $M$ containing $P$, such that $K_{M}L\subseteq
P$, then $K=P$ or $L=P$.

\bigskip

\bigskip The following definition was given in [4, Definition 3.2], and we
include it here for the convenience of the reader.

\bigskip

\textbf{Definition 1.5}. Let $M\in R$- Mod. If $N$ is a submodule of $M$,
then successive powers of $N$ are defined as follows: First, $N^{2}=N_{M}N$.
Then by induction, for any integer $k>2$, we define $N^{k}=N_{M}\left(
N^{k-1}\right) $.

\bigskip

\bigskip Notice that in general is false that $N_{M}\left( N^{k-1}\right)
=\left( N^{k-1}\right) _{M}N$. In the Example 2.8 \ (1) we can see this.

\bigskip

\bigskip

\textbf{Lemma 1.6}. Let $M\in R$-Mod be projective in $\sigma \left[ M\right]
$ and $N\subseteq M$ semiprime in $M$ submodule. If $J$ is a submodule of $M$
such that $J^{n}\subseteq N$ for some positive positive integer $n$, then $%
J\subseteq N$.

\bigskip

\textbf{Proof.} In case $n=1$, there is nothing to prove. Now let $n>1$ and
assume the Lemma holds for lower powers. Since $n\geq 2$, then $2n-2\geq n$.
Hence $J^{2n-2}\subseteq J^{n}$. Now by [1, Proposition 5.6] \ we have that $%
\left( J^{n-1}\right) ^{2}=J^{2n-2}\subseteq N$. Thus by Remark 1.2 we have
that $J^{n-1}\subseteq N$. So by the induction hypothesis we have that $%
J\subseteq N$.

\bigskip

\bigskip

\textbf{Lemma 1.7}. Let $M\in R$-Mod be projective in $\sigma \left[ M\right]
$ and semiprime module . If $N$ and $L$ are submodules of $M$ such that $%
N_{M}L=0$, then

$i)$ \ $L_{M}N=0$

$ii)$ $N\cap L=0$

\bigskip

\textbf{Proof}. $i)$ Since $M$ is projective in $\sigma \left[ M\right] $,
then by [1 Proposition 5.6] we have that $\left( L_{M}N\right) _{M}\left(
L_{M}N\right) =L_{M}\left( N_{M}L\right) _{M}N=0$. As $M$ is semiprime and
by Remark 1.2, then $L_{M}N=0$.

$ii)$ By [3, Proposition 1.3] we have that $\ \ \ \left( N\cap L\right)
_{M}\left( N\cap L\right) \subseteq N_{M}L=0$. Newly as $M$ is semiprime,
then $N\cap L=0$.

\bigskip

\textbf{Remark 1.8.} Let $M\in R$-Mod be projective in $\sigma \left[ M%
\right] $ and a semiprime module, then we claim that $L\cap Ann_{M}\left(
L\right) =0$ for all $L$ submodule of $M$. In fact, if $K=L\cap
Ann_{M}\left( L\right) $ , then by [3, Proposition 1.3] we have that $%
K_{M}K\subseteq Ann_{M}\left( L\right) _{M}L=0$. Since $M$ is semiprime,
then $K=0$. Moreover we claim that $L\cap L^{\prime }\neq 0$ for all $%
L^{\prime }\subseteq _{FI}$ $M$ such that $Ann_{M}\left( L\right)
\varsubsetneq L^{\prime }$. In fact, let $L^{\prime }$ $\subseteq _{FI}M$ be
such that $Ann_{M}\left( L\right) \varsubsetneq L^{\prime }$. Suppose that $%
L\cap L^{\prime }=0$, then $L_{M}^{\prime }L\subseteq L^{\prime }\cap L=0$.
Hence $L^{\prime }\subseteq Ann_{M}\left( L\right) $. So $Ann_{M}\left(
L\right) =L^{\prime }$ a contradiction.

\bigskip

\textbf{Proposition 1.9}. Let $M\in R$-Mod be projective in $\sigma \left[ M%
\right] $ and a semiprime module, then $Hom_{R}\left( M,N\right) \neq 0$ for
all $0\neq N\subseteq M$.

\bigskip

\textbf{Proof}. Suppose that $Hom_{R}\left( M,N\right) =0$. Since $%
Ann_{M}\left( N\right) =$

$\cap \left\{ \ker f\mid f\in Hom_{R}\left( M,N\right) \right\} $, then $%
Ann_{M}\left( N\right) =M$. Hence $M_{M}N=0$. Now by Lemma 1.7 we have that $%
N_{M}M=0$. As $N_{M}M=\sum_{f\in End_{R}\left( M\right) }f\left( M\right) $,
then $N\subseteq N_{M}M$. So $N=0$ a contradiction.

\bigskip

Notice that if $M$ is as in Proposition 1.9, then $\ Ann_{M}\left( N\right)
=M$ if and only if $Hom_{R}\left( M,N\right) =0$ if and only if $N=0$.
Therefore $Ann_{M}\left( N\right) \neq M$ for all $0\neq N\subseteq M$.

\bigskip

\bigskip \textbf{Proposition 1.10.} Let $M\in R$-Mod and $\tau \in M$-tors.
If $P$ is a prime in $M$ and $P$ is $\tau $-pure submodule of $M$, then
there exists $P^{\prime }\subseteq P$ such that $P^{\prime }$ is $\tau $%
-pure and minimal $\ $prime in $M.$

\bigskip

\bigskip

\textbf{Proof}$.$Let $\mathcal{X}=\left\{ Q\in M\mid Q\subseteq P\text{ and }%
Q\text{ is }\tau \text{-pure and prime in }M\text{ \ }\right\} $. Since $%
P\in \mathcal{X}$, then $\mathcal{X}\neq \emptyset $. We claim that any
chain $\mathcal{Y}\subseteq \mathcal{X}$, has a lower bound in $\mathcal{X}$%
. In fact, let $\mathcal{Y=}\left\{ Q_{i}\right\} $ a descending chain chain
in $\mathcal{X}$. Let $Q=\cap \mathcal{Y}$ $\ $is clear that $Q$ is fully
invariant submodule of $M$ and $Q\subseteq P$. \ Since each $Q_{i}$ is $\tau
$-pure in $M$, then $Q$ is $\tau $-pure in $M$. Now $K$ and $L$ be fully
invariant submodule of $M$ such that $K_{M}L\subseteq Q$. If $L\nsubseteq Q$%
, then there exists $Q_{t}\in \mathcal{Y}$ such that $L\nsubseteq Q_{t}$.
Hence $L\nsubseteq Q_{l}$ for all $Q_{l}\subseteq Q_{t}$. Thus $K\subseteq
Q_{l}$ for all $Q_{l}\subseteq Q_{t}$. As $\mathcal{Y}$ is a descending
chain, then $K\subseteq Q$. Therefore by Zorn's Lemma $\mathcal{X}$ has
minimal elements.

\bigskip

\bigskip

\textbf{Proposition 1.11.} Let $M$ $\in R$-Mod be projective in $\sigma %
\left[ M\right] $, semiprime module and $\tau \in M$-tors. If $U\subseteq M$
is a uniform submodule, such that $U\in \mathbb{F}_{\tau }$, then $%
Ann_{M}\left( U\right) $ is $\tau $-pure and prime in $M$.

\bigskip

\textbf{Proof}$.$ By Proposition 1.9 we have that $Ann_{M}\left( U\right)
\varsubsetneq M$. Now we claim that $Ann_{M}\left( U\right) =Ann_{M}\left(
U^{\prime }\right) $ for all $0\neq U^{\prime }\subseteq U$. In fact let $%
0\neq U^{\prime }\subseteq U$, then $Ann_{M}\left( U\right) \subseteq
Ann_{M}\left( U^{\prime }\right) $. Suppose that $Ann_{M}\left( U\right)
\varsubsetneq Ann_{M}\left( U^{\prime }\right) $. By [3, Proposition 1.9] we
have that $Ann_{M}\left( U^{\prime }\right) $ is a fully invariant submodule
of $M$. Now by Remark 1.8 we have that $U\cap Ann_{M}\left( U^{\prime
}\right) \neq 0$. As $U$ is uniform module, then $\left[ U\cap Ann_{M}\left(
U^{\prime }\right) \right] \cap U^{\prime }\neq 0$. But by Remark 1.8 we
have that $Ann_{M}\left( U^{\prime }\right) \cap U^{\prime }=0$. So $\left[
U\cap Ann_{M}\left( U^{\prime }\right) \right] \cap U^{\prime }=0$ a
contradiction. Therefore $Ann_{M}\left( U\right) =Ann_{M}\left( U^{\prime
}\right) $. \ Now by [3, Lemma 1.16] $Ann_{M}\left( U\right) $ is prime in $%
M $. On the other hand we know that $Ann_{M}\left( U\right) =\cap \left\{
\ker f\mid f\in Hom_{R}\left( M,U\right) \right\} $. Since $U\in \mathbb{F}%
_{\tau }$, $\ker f$ is $\tau $-pure in $M$. Hence $Ann_{M}\left( U\right) $
is $\tau $-pure in $M$.

\bigskip

\bigskip

\textbf{Remark 1.12.} \ Let $M$ be as in the Proposition 1.11 and $%
U\subseteq M$ \ uniform submodule, then $Ann_{M}\left( U\right)
=Ann_{M}\left( U^{\prime }\right) =P$ for all $0\neq U^{\prime }\subseteq U$%
, where $P$ is prime in $M$. So by [3, Definition 4.3] we have that $%
Ass_{M}\left( U\right) =\{P\}$. \ Moreover by [3, Proposition 4.4] we have
that $Ass_{M}\left( U^{\prime }\right) =Ass_{M}\left( U\right) =\left\{
P\right\} $ for all $0\neq U^{\prime }\subseteq U$.

\bigskip

\bigskip

\textbf{Lemma 1.13}. Let $M\in R$-Mod be projective in $\sigma \left[ M%
\right] $ and semiprime module. Suppose that \ $N$ is a submodule of $M$
such that $P=Ann_{M}\left( N\right) $ is prime in $M$, then $P$ is minimal
prime in $M$.

\bigskip

\textbf{Proof.} By Proposition 1.10, there exists $P^{\prime }$ a minimal
prime in $M$ such that $P^{\prime }\subseteq P$. As $\left[ Ann_{M}\left(
N\right) \right] _{M}N=0\subseteq P^{\prime }$, then $N\subseteq P^{\prime }
$ or $Ann_{M}\left( N\right) \subseteq P^{\prime }$. If $N\subseteq
P^{\prime }$, then $N\subseteq P=Ann_{M}\left( N\right) $. Thus $N_{M}N=0$.
Since $M$ is semiprime module, then $N=0$. Hence $P=Ann_{M}\left( N\right)
=Ann_{M}\left( 0\right) =M$, it is not possible. So the only possibility is $%
P=Ann_{M}\left( N\right) \subseteq P^{\prime }$. Therefore $P=P^{\prime }$
and we have the result.

\bigskip

Note that if $M$ is as in the Lemma 1.13 and $U$ is a uniform submodule of $%
M $, then by Proposition 1.11 we have that $P=Ann_{M}\left( U\right) $ is a
prime in $M$ submodule. Now by Lemma 1.13 we obtain that $P=Ann_{M}\left(
U\right) $ is minimal prime in $M$. Moreover by proof of the Proposition
1.11 we have that $P=Ann_{M}\left( U\right) =Ann_{M}\left( \widehat{U}%
\right) $.

\bigskip

\bigskip

\textbf{Definition 1.14}. Let $M\in R$-Mod be projective in $\sigma \left[ M%
\right] $. If $N$ is a submodule of $M$, the right annihilator of $N$ in $M$
is $Ann_{M}^{r}\left( N\right) =\sum \left\{ L\subseteq M\mid
N_{M}L=0\right\} $.

\bigskip

Note that $Ann_{M}^{r}\left( N\right) $ is the largest submodule of $M$ such
that

$N_{M}\left[ Ann_{M}^{r}\left( N\right) \right] =0$. In fact, if $L$ and $K$
are submodules of $M$ such that $N_{M}L=0$ and $N_{M}K=0$, then by [3,
Proposition 1.3] we have that $N_{M}\left( L\oplus K\right) =0$. Now since $%
M $ is projective in $\sigma \left[ M\right] $ then by [1, Proposition 5.5] $%
N_{M}\left( \dfrac{L\oplus K}{T}\right) =0$ for all $T\subseteq L\oplus K$.
So $N_{M}\left( L+K\right) =0$.

\bigskip

\textbf{Remark 1.15}. If $M$ is an $R$-module projective\ in $\sigma \left[ M%
\right] $ and $N$ and $L$ are submodules of $M$ such that $N_{M}L=0$, then
we claim that $N_{M}\left( L_{M}M\right) =0$. In fact by [1, Proposition
5.6] we have that $N_{M}\left( L_{M}M\right) =\left( N_{M}L\right) _{M}M=0$.
Since $L\subseteq L_{M}M$ and $L_{M}M$ is a fully invariant submodule of $M$%
, then $Ann_{M}^{r}\left( N\right) =\sum \left\{ L\subseteq _{FI}M\mid
N_{M}L=0\right\} $. Hence $Ann_{M}^{r}\left( N\right) $ is a fully invariant
submodule of $M$.

\bigskip

\textbf{Proposition 1.16}. Let $M\in R$-Mod be projective in $\sigma \left[ M%
\right] $. If $M$ is a semiprime module, then $Ann_{M}\left( N\right)
=Ann_{M}^{r}\left( N\right) $ for all $N\subseteq M$.

\bigskip

\textbf{Proof}. Since $Ann_{M}\left( N\right) _{M}N=0$, then by Lemma 1.7, $%
N_{M}Ann_{M}\left( N\right) =0$. So $Ann_{M}\left( N\right) \subseteq
Ann_{M}^{r}\left( N\right) $. Newly as $N_{M}Ann_{M}^{r}\left( N\right) =0$,
then

$Ann_{M}^{r}\left( N\right) _{M}N=0$. Thus $Ann_{M}^{r}\left( N\right)
\subseteq Ann_{M}\left( N\right) $ Hence $Ann_{M}\left( N\right)
=Ann_{M}^{r}\left( N\right) $.

\bigskip

\bigskip

\textbf{Definition} \textbf{1.17}. Let $M\in R$-Mod. A submodule $N$ of $M$
is named annihilator submodule, if $N=Ann_{M}\left( K\right) $ for some $%
K\subseteq M$.

\bigskip

Notice that if $N$ is an annihilator submodule, then $N$ is a fully
invariant submodule of $M$.

\bigskip

\textbf{Proposition 1.18}. Let $M\in R$-Mod be projective in $\sigma \left[ M%
\right] $. If $M$ is a semiprime module, then $Ann_{M}\left( Ann_{M}\left(
N\right) \right) =Ann_{M}^{r}\left( Ann_{M}^{r}\left( N\right) \right) =N$
for all $N$ annihilator submodule of $M$.

\bigskip

\textbf{Proof. }As $\ N_{M}Ann_{M}^{r}\left( N\right) =0$, then by Lemma 1.7
we have that

$Ann_{M}^{r}\left( N\right) _{M}N=0$. Thus $N\subseteq Ann_{M}^{r}\left(
Ann_{M}^{r}\left( N\right) \right) $. Since $N=Ann\left( K\right) $ for some
submodule of $M$, then $N_{M}K=0$. Hence $K\subseteq Ann_{M}^{r}\left(
N\right) $. Therefore $Ann_{M}^{r}\left( Ann_{M}^{r}\left( N\right) \right)
\subseteq Ann_{M}^{r}\left( K\right) $. So by Proposition 1.16 we have that

$Ann_{M}^{r}\left( Ann_{M}^{r}\left( N\right) \right) \subseteq
Ann_{M}\left( K\right) =N$. Hence $Ann_{M}^{r}\left( Ann_{M}^{r}\left(
N\right) \right) =N$. Moreover by Proposition 1.16 we have that

$Ann_{M}\left( Ann_{M}\left( N\right) \right) =Ann_{M}^{r}\left(
Ann_{M}^{r}\left( N\right) \right) $.

\bigskip

\section{\protect\bigskip Modules with ACC on annihilators and Goldie Modules%
}

\bigskip

\bigskip The following definition was given in [4, Definition 3.1]. We
include here this definition for the convenience of the reader.

Let $M\in R$-Mod. For a subset $X\subseteq End_{R}\left( M\right) $, let $%
\mathcal{A}_{X}=\cap \left\{ \ker f\mid \text{ }f\in X\right\} $. Now we
consider the set $\mathcal{A}_{M}=\left\{ \mathcal{A}_{X}\mid \text{ }%
X\subseteq End_{R}\left( M\right) \text{ }\right\} $.

\bigskip

\textbf{Definition 2.1}. Let $M\in R$-Mod, we say $M$ satisfies ascending
chain condition (ACC) on annihilators , if \ $\mathcal{A}_{M}$ satisfies ACC.

\bigskip

Notice that if $M=R$, then $R$ satisfies ACC on annihilators in the sense of
Definition 2.1, \ if and only if $R$ satisfies ACC on left annihilators in
the usual sense for a ring $R$.

\bigskip

Also note that if $K$ is a submodule of $M$, we have that

$Ann_{M}\left( K\right) =\underset{f\in X}{\cap }\left\{ \ker f\mid
X=Hom_{R}\left( M,K\right) \right\} $. \

Since $Hom_{R}\left( M,K\right) \subseteq End_{R}\left( M\right) $, then $%
Ann_{M}\left( K\right) =\mathcal{A}_{X}$,

where $X=Hom_{R}\left( M,K\right) $. Hence $Ann_{M}\left( K\right) \in
\mathcal{A}_{M}$

\bigskip

\textbf{Theorem 2.2.} \ Let $M\in R$-Mod be projective in $\sigma \left[ M%
\right] $ and a semiprime module. If $M$ satisfies ACC on annihilators, then:

\bigskip

$i)$ $M$ has finitely many minimal of minimal prime in $M$ submodules.

$ii)$ If $P_{1}$, $P_{2}$, ..., $P_{n}$ are the minimal prime in $M$
submodules, then $P_{1}\cap $ $P_{2}\cap $ , ..., $\cap P_{n}=0$

$iii)$ If $P\subseteq M$ is prime in $M$, then $P$ is minimal prime in $M$
if and only if $P$ is an annihilator submodule.

\bigskip

\textbf{Proof. \ }$i)$ By a "prime annihilator " we mean a prime in $M$
submodule which is an annihilator submodule. We shall first show that every
annihilator submodule of $M$ contains a finite product of annihilators
primes in $M$. Suppose not, then by hypothesis there is an annihilator
submodule $N$ which is maximal with respect to not containing a finite
product of annihilators primes in $M$. Hence $N$ is not prime in $M$. Thus
by Remark 1.4 there are fully invariant submodules $L$ and $K$ of $M$ such
that $N\varsubsetneq L$, $N\varsubsetneq K$ and $L_{M}K\subseteq N$. Since $%
Ann_{M}\left( N\right) _{M}N=0$, then by [3, Proposition 1.3] $Ann_{M}\left(
N\right) _{M}\left( L_{M}K\right) =0$. Since $M$ is projective in $\sigma %
\left[ M\right] $, then by [ 1, Proposition 5.6] we have that $\left(
Ann_{M}\left( N\right) _{M}L\right) _{M}K=0$. Therefore $K\subseteq
Ann_{M}^{r}\left( Ann_{M}\left( N\right) _{M}L\right) $. Now we claim that $%
L_{M}\left[ Ann_{M}^{r}\left( Ann_{M}\left( N\right) _{M}L\right) \right]
\subseteq N$. We have that

$\left[ Ann_{M}\left( N\right) _{M}L\right] _{M}\left[ Ann_{M}^{r}\left(
Ann_{M}\left( N\right) _{M}L\right) \right] =0$. Newly as $M$ is projective
in $\sigma \left[ M\right] $, then $Ann_{M}\left( N\right) _{M}\left( L_{M}%
\left[ Ann_{M}^{r}\left( Ann_{M}\left( N\right) _{M}L\right) \right] \right)
=0$. Hence

$L_{M}\left[ Ann_{M}^{r}\left( Ann_{M}\left( N\right) _{M}L\right) \right]
\subseteq Ann_{M}^{r}\left( Ann_{M}\left( N\right) \right) $. Since $N$ is
annihilator submodule, then by Lemma 1.7 and Proposition 1.18 we have that

$L_{M}\left[ Ann_{M}^{r}\left( Ann_{M}\left( N\right) _{M}L\right) \right]
\subseteq N$. Since $K\subseteq Ann_{M}^{r}\left( Ann_{M}\left( N\right)
_{M}L\right) =Ann_{M}\left( Ann_{M}\left( N\right) _{M}L\right) $, then we
can take $K=Ann_{M}\left( Ann_{M}\left( N\right) _{M}L\right) $. Thus $K$ is
an annihilator submodule. \ Similarly we can prove that

$L\subseteq Ann_{M}\left[ K_{M}Ann_{M}^{r}\left( N\right) \right] $ and $%
\left( Ann_{M}\left[ K_{M}Ann_{M}^{r}\left( N\right) \right] \right)
_{M}K\subseteq N$. So also we can take

$L=Ann_{M}\left[ K_{M}Ann_{M}^{r}\left( N\right) \right] $. Hence $L$ is an
annihilator submodule. Therefore $K$ and $L$ are annihilator submodules such
that $N\varsubsetneq K$, $N\varsubsetneq L$, and $L_{M}K\subseteq N$. So $L~$%
and $K$ contain a finite product of prime annihilators. Hence $N$ contains a
finite product of prime annihilator, this is a contradiction. Therefore
every annihilator submodule of $M$ contains a finite product of prime
annihilators. Since $0=Ann_{M}\left( M\right) $, then the zero is an
annihilator submodule. Thus there are prime annihilators $P_{1}$, $P_{2}$,
..., $P_{n}$ of $M$ such that $\left( P_{1}\right) _{M}$ $\left(
P_{2}\right) _{M}$, ..., $P_{n}=0$. Now if $Q$ is a minimal prime in $M$,
then $\left( P_{1}\right) _{M}$ $\left( P_{2}\right) _{M}$ ... $_{M}\left(
P_{n}\right) =0\subseteq Q$. Thus there exists $P_{j}\subseteq Q$ for some $%
0\leq j\leq n$. Whence $P_{j}=Q$. Therefore the minimal primes in $M$ are
contained in the finite set $\left\{ P_{1}\text{, }P_{2}\text{, ... , }%
P_{n}\right\} $. Moreover each $P_{i}$ is an annihilator submodule for $%
0\leq i\leq n$.

\bigskip

$ii)$ By $i)$ we can suppose that $P_{1}$, $P_{2}$, ..., $P_{n}$ are the
minimal prime in $M$ submodules. Now by [3 Proposition 1.3] we have that \ \
$\left[ P_{1}\cap P_{2}\cap ,...,\cap P_{n}\right] ^{n}\subseteq \left(
P_{1}\right) _{M}\left( P_{2}\right) _{M}..._{M}\left( P_{n}\right) =0$.
Since $M$ is semiprime module, then by Lemma 1.6 we have that $P_{1}\cap
P_{2}\cap ,...,\cap P_{n}=0$

\bigskip

$iii)$ Let $P\subseteq M$ be a prime in $M$ submodule. If $P$ is minimal
prime in $M$, then by $ii)$ we have that $P=P_{j}$ for some $1\leq j\leq n$.
By proof of $i)$ we know that each $P_{t}$ is an annihilator submodule for $%
1\leq t\leq n$. Hence $P$ so does. Now we obtain the converse by By Lemma
1.13.

.

\bigskip

Notice that in Theorem 2.2$~$ $i)$ we only use the condition $M$ satisfies
ACC on annihilator submodules. Also note that in $iii)$ each minimal prime $%
P_{i}$ in $M$ is an annihilator submodule. Thus $P_{i}=Ann_{M}\left(
K\right) $ for some $K\subseteq \dot{M}$. As $K\in \mathbb{F}_{\chi \left(
M\right) }$ and $Ann_{M}\left( K\right) =\cap \left\{ \ker f\mid f\in
Hom_{R}\left( M,K\right) \right\} $, then $P_{i}$ is $\chi \left( M\right) $%
-pure in $M$ for all $1\leq i\leq n$.

\bigskip

Note that in Theorem 2.2 $M$ satisfies ACC on annihilators is a necessary
condition. In order to see this, consider de following example

\bigskip

\textbf{Example 2.3.} Let $R$ a ring such that $\left\{ S_{i}\right\} _{i\in
I}$ is a family of non isomorphic simple $R$-modules and let $%
M=\bigoplus_{i\in I}S_{i}$ . It is clear that $M$ is projective in $\sigma %
\left[ M\right] $. Now let $N\subseteq M$ such that $N_{M}N=0$. Since $M$ is
semisimple module, then $N=0$. Thus $M$ is a semiprime module. We claim that
$M$ does not satisfy ACC on annihilators. In fact, let $L$ be a submodule of
$M$, then there exists $K$ such that $M=K\oplus L$. So we can consider the\
canonical projection $\pi :M\rightarrow K$. Hence $\ker \pi =L$. Thus each
submodule of $M$ is the kernel of some morphism. Hence $M$ does not satisfy
ACC on annihilators. Now since $S_{i}\ncong S_{j}$, then $S_{i}$ is minimal
prime in $M$. So $M$ has infinitely many minimal prime in $M$ submodules.

\bigskip

\textbf{Lemma 2.4}. Let $M\in R$-Mod. If $C\in \sigma \left[ M\right] $ is $%
\chi \left( M\right) $-cocritical, then there are submodules $C^{\prime
}\subseteq C$ and $M^{\prime }\subseteq M$ such that $C^{\prime }$ is
isomorphic to $M^{\prime }$.

\bigskip

\textbf{Proof}. Since $C$ is $\chi \left( M\right) $-cocritical, then $C\in
\mathbb{F}_{\chi \left( M\right) }$. Thus $Hom_{R}\left( C,\text{ }\widehat{M%
}\right) \neq 0$. So there exists $\ $ $f:C\rightarrow $ $\widehat{M}$\ \ \
a non zero morphism. Since $C$ is $\chi \left( M\right) $-cocritical and $%
M\subseteq _{ess}\widehat{M}$, then there exists $C^{\prime }$ submodule of $%
C$ such that $C^{\prime }\hookrightarrow M$. Hence there exists $M^{\prime
}\subseteq M$ such that $C^{\prime }\cong M^{\prime }$. Also note that $%
M^{\prime }$ is $\chi \left( M\right) $-cocritical.

\bigskip

\bigskip

\textbf{Remark 2.5}. If $M$ is as in the Theorem 2.2 and $C\in \sigma \left[
M\right] $ is $\chi \left( M\right) $-cocritical, then by Lemma 2.4 there
are $C^{\prime }\subseteq C$ and $M^{\prime }\subseteq M$ such that $%
C^{\prime }\cong M^{\prime }$. Hence $M^{\prime }$ is a uniform module, then
by Proposition 1.11 we have that $Ann_{M}\left( M^{\prime }\right) =P$ where
$P$ is prime in $M$. By Lemma 2.13 $P$ is minimal prime in $M$. Moreover by
Remark 1.12 we have that $Ass_{M}\left( M^{\prime }\right) =\left\{
P\right\} $. Hence $Ass_{M}\left( C^{\prime }\right) =\{P\}$. Since $%
C^{\prime }\subseteq _{ess}C$, then by [3, Proposition 4.4] we obtain $%
Ass_{M}\left( C\right) =\{P\}$.

\bigskip

\bigskip

Let $Spec^{Min}\left( M\right) $ denote the set of minimal primes in $M$ and
$\mathcal{E}_{\tau }\left( M\right) $ a complete set of representatives of
isomorphism classes of indecomposable $\tau $-torsion free injective modules
in $\sigma \left[ M\right] $.

Also we denote by $\mathcal{P}_{M}^{Min}=\left\{ \chi \left( M/P\right) \mid
P\subseteq M\text{ is minimal prime in }M\right\} $.

\bigskip

\textbf{Remark 2.6}. Let $\tau _{g}$ be the hereditary torsion theory
generated by the family of $M$-singular modules in $\sigma \left[ M\right] $%
.\ If $M$ is non $M$-singular, then $\tau _{g}=\chi \left( M\right) $. In
fact, if $M$ is non $M$-singular, then $\tau _{g}\leq \chi \left( M\right) $%
. Now if $\tau _{g}<\chi \left( M\right) $, then there exists $0\neq N\in
\mathbb{T}_{\chi \left( M\right) }$ such that $N\in \mathbb{F}_{\tau _{g}}$.
Thus $Hom_{R}\left( N,\widehat{M}\right) =0$. By [17, Proposition 10.2] $N$
is $M$-singular. Thus $N\in $ $\mathbb{T}_{\tau _{g}}$ a contradiction. So $%
\tau _{g}=\chi \left( M\right) $.

\bigskip

\textbf{Theorem 2.7 }. Let $M\in R$-Mod be projective in $\sigma \left[ M%
\right] $ and semiprime module. If $M$ satisfies ACC on annihilators and for
each $0\neq N\subseteq M$, $N$ contains a uniform submodule, then\ there is
a bijective correspondences between $\mathcal{E}_{\chi \left( M\right)
}\left( M\right) \,$and $Spec^{Min}\left( M\right) $.

\bigskip

\textbf{Proof}. Since $M$ is semiprime and $M$ satisfies ACC on annihilators
then by [4, Proposition 3.4] we have that $M$ $\ $is non $M$-singular. Thus
by Remark 2.6 $\tau _{g}=\chi \left( M\right) $. Let $E\in \mathcal{E}%
_{_{\chi \left( M\right) }}\left( M\right) $. As $E$ is a uniform module and
$E\in \mathbb{F}_{\chi \left( M\right) }$, then $E$ is $\chi \left( M\right)
$-cocritical. Now by Remark \textbf{\ }2.5, we have that $Ass_{M}\left(
E\right) =\left\{ P\right\} $ with $P\in Spec^{Min}\left( M\right) $. Hence
we define the function

\begin{equation*}
\Psi :\mathcal{E}_{\chi \left( M\right) }\left( M\right) \rightarrow
Spec^{Min}\left( M\right)
\end{equation*}%
as $\Psi \left( E\right) =P$. We claim that $\Psi $ is bijective. Suppose
that $\Psi \left( E_{1}\right) =\Psi \left( E_{2}\right) =P$. \ Since $E_{1}$
and $E_{2}$ are uniform modules, then $E_{1}$ and $E_{2}$ are $\chi \left(
M\right) $-cocritical. \ Since $\tau _{g}=\chi \left( M\right) $, then by
Lemma 2.4 there are $C_{1}^{\prime }\subseteq E_{1}$, $C_{2}^{\prime
}\subseteq E_{2}$ and $M_{1}$, $M_{2}$ submodules of $M$ such that $%
C_{1}^{\prime }\cong M_{1}$ and $C_{2}^{\prime }\cong M_{2}$. Hence $M_{1}$
and $M_{2}$ are uniform modules, then by Remark 1.12 we have that. $%
Ann_{M}\left( C_{1}^{\prime }\right) =Ann_{M}\left( M_{1}\right)
=P=Ann_{M}\left( M_{2}\right) =Ann_{M}\left( C_{2}^{\prime }\right) $. So $%
Ann_{M}\left( M_{1}\right) =Ann_{M}\left( M_{2}\right) =P$.

On the other hand we have that $\left( M_{1}+P\right) /P$ and $\left(
M_{2}+P\right) /P$ are submodules of $M/P$. \ By [3, Proposition 2.2 and 2.7
] $\chi \left( \left( M_{1}+P\right) /P\right) =\chi \left( M/P\right) =\chi
\left( \left( M_{2}+P\right) /P\right) $. Since $Ann_{M}\left( M_{1}\right)
=Ann_{M}\left( M_{2}\right) =P$, then by Remark 1.8, $M_{1}\cap P=0$ \ and $%
M_{2}\cap P=0$. Therefore $\left( M_{1}+P\right) /P\cong M_{1}/\left( P\cap
M_{1}\right) =M_{1}$ and $\left( M_{2}+P\right) /P\cong M_{2}/\left( P\cap
M_{1}\right) =M_{2}$. Thus $\chi \left( M_{1}\right) =\chi \left( M/P\right)
=\chi \left( M_{2}\right) $. So $Hom_{R}\left( M_{1},\widehat{M_{2}}\right)
\neq 0$. Since $M_{1}$ and $M_{2}$ are $\tau _{g}$-cocritical, then there
exists $N_{1}\subseteq M_{1}$ such that $N_{1}\hookrightarrow M_{2}$. As $%
M_{1}$ and $M_{2}$ are uniform modules, then $\widehat{M_{1}}\cong $ $%
\widehat{N_{1}}\cong $ $\widehat{M_{2}}$. Thus $E_{1}=\widehat{C_{1}^{\prime
}}$ $\cong \widehat{M_{1}}\cong $ $\widehat{M_{2}}\cong \widehat{%
C_{2}^{\prime }}=E_{2}$. So $\Psi $ is injective.

\bigskip

Now let $P\in Spec^{Min}\left( M\right) $ since $M$ satisfies ACC on
annihilators, then by Theorem 2.2 $P=Ann_{M}\left( N\right) $ for some $%
0\neq N\subseteq M$. By hypothesis there exists a uniform module $U$ such
that $U\subseteq N$. Thus $P=Ann_{M}\left( N\right) \subseteq Ann_{M}\left(
U\right) $. By Proposition 1.11, we have that $Ann_{M}\left( U\right) $ is
prime in $M$. As $Ann_{M}\left( U\right) $ is clearly an annihilator
submodule, then by Theorem 2.2 $iii)$ we have that $Ann_{M}\left( U\right) $
is minimal prime in $M$. Therefore $P=Ann_{M}\left( U\right) $. Furthermore
by Remark 1.12 we have that $Ass_{M}\left( U\right) =\{P\}$. Hence $%
Ass_{M}\left( \widehat{U}\right) =\{P\}$ Thus $\Psi \left( \widehat{U}%
\right) =P$. As $U\subseteq M$, then $U\in F_{\chi \left( M\right) }$.
Therefore $\widehat{U}$ is an indecomposable $\chi \left( M\right) $-torsion
free injective module in $\sigma \left[ M\right] $. Hence $\Psi $ is
surjective.

\bigskip

Notice that if $P$ is prime in $M$ such that $P$ is $\chi \left( M\right) $%
-pure, then $P$ is minimal prime in $M$. In fact by Proposition 1.10 there
exists $P^{\prime }\subseteq P$ such that $P^{\prime }$ is $\chi \left(
M\right) $-pure and minimal prime in $M$. Now let $N$ be a submodule of $M$
such that $P\cap N=0$, then $P_{M}N\subseteq P\cap N=0$. So $P_{M}N\subseteq
P^{\prime }$. Thus $P\subseteq P^{\prime }$or $N\subseteq P^{\prime }$.
Suppose that $P\nsubseteq P^{\prime }$, then $N\subseteq P^{\prime }$. Thus $%
N\subseteq P$. Hence $0=P\cap N=N$. So $P\subseteq _{ess}M$. Thus $M/P\in
\mathbb{T}_{\chi \left( M\right) }$ a contradiction. So $P\subseteq
P^{\prime }$. Thus $P=P^{\prime }$. On the other hand , if $Spec_{\chi
\left( M\right) }\left( M\right) $ denotes the set of $\chi \left( M\right) $%
-pure submodules prime in $M$, then $Spec^{Min}\left( M\right) =Spec_{\chi
\left( M\right) }\left( M\right) $. So if $M$ is as in the Theorem 2.7, then
$M$ has local Gabriel correspondence with respect to $\chi \left( M\right) $.

\bigskip

\bigskip \textbf{Example 2.8. }Let $R=\mathbb{Z}_{2}\rtimes \left( \mathbb{Z}%
_{2}\oplus \mathbb{Z}_{2}\right) $, the trivial extension of $\mathbb{Z}_{2}$
by $\mathbb{Z}_{2}\oplus \mathbb{Z}_{2}$. This ring can be described as

$R=\left \{ \left(
\begin{array}{cc}
a & \left( x,y\right) \\
0 & a%
\end{array}%
\right) \mid a\in \mathbb{Z}_{2},\left( x,y\right) \in \mathbb{Z}_{2}\oplus
\mathbb{Z}_{2}\right \} $. $R$ has only one maximal ideal $I=\left \{ \left(
\begin{array}{cc}
0 & \left( x,y\right) \\
0 & 0%
\end{array}%
\right) \mid \left( x,y\right) \in \mathbb{Z}_{2}\oplus \mathbb{Z}%
_{2}\right
\} $ and it has three simple ideals: $J_{1}$, $J_{2}$, $J_{3}$
which are isomorphic, where $J_{1}=\left \{ \left(
\begin{array}{cc}
0 & \left( 0,0\right) \\
0 & 0%
\end{array}%
\right) ,\left(
\begin{array}{cc}
0 & \left( 1,0\right) \\
0 & 0%
\end{array}%
\right) \right \} $, $J_{2}=\left \{ \left(
\begin{array}{cc}
0 & \left( 0,0\right) \\
0 & 0%
\end{array}%
\right) ,\left(
\begin{array}{cc}
0 & \left( 0,1\right) \\
0 & 0%
\end{array}%
\right) \right \} $, $J_{3}=\left \{ \left(
\begin{array}{cc}
0 & \left( 0,0\right) \\
0 & 0%
\end{array}%
\right) ,\left(
\begin{array}{cc}
0 & \left( 1,1\right) \\
0 & 0%
\end{array}%
\right) \right \} $. Then the lattice of ideals of $R$ has the following form

\bigskip

\begin{equation*}
\begin{array}{ccc}
& \overset{R}{\bullet } &  \\
& \overset{I}{\bullet } &  \\
\overset{J_{1}}{\bullet } & \overset{J_{2}}{\bullet } & \overset{J_{3}}{%
\bullet } \\
& \overset{0}{\bullet } &
\end{array}%
\end{equation*}

\bigskip

$R$ is artinian and $R$-$Mod$ \ has only one simple module up to
isomorphism. Let $S$ be the simple module. By [12, Theorem 2.13], \ we know
that there is a lattice anti-isomorphism between the lattice of ideals of $R$
and the lattice of fully invariant submodules of $E\left( S\right) $. Thus
the lattice of fully invariant submodules of $E\left( S\right) $ has tree
maximal elements $K$, $L$ and $N$. Moreover, $E\left( S\right) $ contains
only one simple module $S$. Therefore, the lattice of fully invariant
submodules of $E\left( S\right) $ has the following form.

\bigskip

\begin{equation*}
\begin{array}{ccc}
& \overset{E\left( S\right) }{\bullet } &  \\
\overset{K}{\bullet } & \overset{L}{\bullet } & \overset{N}{\bullet } \\
& \overset{S}{\bullet } &  \\
& \overset{0}{\bullet } &
\end{array}%
\end{equation*}

\bigskip

Let $M=E\left( S\right) $. As $K\cap L=S$ and $K_{M}L\subseteq K\cap L$,
then $K_{M}L\subseteq S$. On the other hand, we consider the morphism. $f:M%
\overset{\pi }{\rightarrow }\left( M/N\right) \cong S$ $\overset{i}{%
\hookrightarrow }L$, where $\pi $ is the canonical projection \ and $i$ is
the inclusion. Thus $f\left( K\right) =S$. Therefore $S\subseteq K_{M}L$.
Thus we have that $K_{M}L=S$. Thus $K_{M}L\subseteq N$ but $K\nsubseteq N$
and $L\nsubseteq N$. Therefore $N$ is not prime in $M$. Analogously, we
prove that neither $K$ nor $L$ are prime in $M$. Also we note that $K_{M}K=S$
and $S_{M}K=0$. Thus $A_{nn_{M}}\left( K\right) =S$.

We claim that $K_{M}S=S$. In fact, if $\pi $ is the natural projection from $%
M$ onto $M/N$, then $\pi \left( K\right) =S$ since $\left( M/N\right) \cong
S $. Thus $K_{M}S=S$. Analogously we prove that $L_{M}S=S$ and $N_{M}S=S$.
Moreover, it is not difficult to prove that $A_{nn_{M}}\left( S\right)
=K\cap L\cap N=S$ and that $S$ is not a prime in $M$ submodule. $Soc\left(
R\right) =J_{1}+$ $J_{2}+J_{3}$ $=J_{1}\oplus $ $J_{2}$ and $Soc\left(
R\right) \subset _{ess}R$. As $S$ is the only simple module, then $Soc\left(
R\right) \cong S\oplus S$. So $E\left( S\oplus S\right) =E\left( R\right) $.
Thus $E\left( R\right) =E\left( S\right) \oplus E\left( S\right) =M\oplus M$%
. Therefore $\sigma \left[ M\right] =\sigma \left[ R\right] =R$-$Mod$.
Moreover $Hom_{R}\left( M,K\right) \neq 0$ for all $K\in \sigma \left[ M%
\right] $ but $M$ is not a generator of $\sigma \left[ M\right] $.

Now let $0\neq M^{\prime }$ a submodule of $M$. As $S\subseteq _{ess}M$,
then $S\subseteq M^{\prime }$. Therefore every submodule of $M$ contains a
uniform submodule of $M$.

From the previous arguments we can conclude that:

1) $K_{M}\left( K_{M}K\right) =K_{M}S=S$ and $\left( K_{M}K\right)
_{M}K=S_{M}K=0$.

Thus $K_{M}\left( K_{M}K\right) \neq \left( K_{M}K\right) _{M}K$

2) $M$ does not have prime in $M$ submodules. Thus $Spec^{Min}\left(
M\right) =\emptyset $

3) $M$ is the only one indecomposable injective module in $\sigma \lbrack M]$%
. Thus $\mathcal{E}_{\chi \left( M\right) }\left( M\right) =\left\{
M\right\} $

4) $M$ is noetherian. So $M$ satisfies ACC on annihilators

5) $M$ is not projective in $\sigma \left[ M\right] $

\bigskip

\textbf{Remark 2.9}. If $M$ is projective in $\sigma \lbrack M]$ and we
define the function $\Phi :Spec^{Min}\left( M\right) \rightarrow $ $\mathcal{%
P}_{M}^{Min}$ as $\Phi \left( P\right) =\chi \left( M/P\right) $, then $\Phi
$ is bijective. In fact. Let $P$ and $P^{\prime }$ are minimal primes in $M$
submodules, such that $\chi \left( M/P\right) =\chi \left( M/P^{\prime
}\right) $, then $P^{\prime }$ is not $\chi \left( M/P\right) $-dense in $M$%
. As $P^{\prime }$ is fully invariant submodule of $M$, then by [3, Lemma
2.6] we obtain $P^{\prime }\subseteq P$. But $P$ is minimal prime in $M$,
then $P^{\prime }=P$. So $\Phi $ is injective. Moreover it is clear that $%
\Phi $ is surjective.

\bigskip

\textbf{Corollary 2.10. }Let $M\in R$-Mod be projective in $\sigma \left[ M%
\right] $ and semiprime module. If $M$ satisfies ACC on annihilators and for
each $0\neq N\subseteq M$, $N$ contains a uniform submodule, then\ there is
a bijective correspondence between $\mathcal{E}_{\chi \left( M\right)
}\left( M\right) \,$and $\mathcal{P}_{M}^{Min}$.

\bigskip

\textbf{Proof.} It follows\ from Theorem 2.7 and Remark 2.9.

\bigskip

\bigskip

\textbf{Definition 2.11}. An $R$-module $M$ is Goldie Module if it satisfies
ACC on annihilators and it has finite uniform dimension.

\bigskip

Notice that if $M=R$, then $R$ is Goldie module in the sense of the
definition 2.11 if and only if $R$ is left Goldie ring$.$

\bigskip

\textbf{Corollary 2.12}. Let $M\in R$-Mod be projective in $\sigma \left[ M%
\right] $ and a semiprime module. If $M$ is a Goldie module, then\ there is
a bijective correspondences between $\mathcal{E}_{\tau _{g}}\left( M\right)
\,$and $Spec^{Min}\left( M\right) $.

\bigskip

\textbf{Proof}. It follows from Definition 2.11 and Theorem\textbf{\ }2.7.

\bigskip

\textbf{Corollary 2.13}. Let $R$ be a semiprime ring such that satisfies ACC
on left annihilators. Suppose that for each $0\neq I\subseteq R$ left ideal
of $R$, $I$ contains a uniform left ideal of $R$, then there is a bijective
correspondence between the set of representatives of isomorphism classes of
indecomposable non singular injective $R$-modules and the set of minimal
prime ideals of $R$.

\bigskip

\textbf{Proof}. From [3, Definition 1.10] we have that $P$ is prime in $R$
if and only if $P$ is prime ideal of $R$. Since $R$ is a semiprime ring and
satisfies ACC on left annihilators, $R$ is non singular

On the other hand we know that $\sigma \left[ R\right] =R$-Mod. So by
Theorem 2.7 we have the result.

\bigskip

Notice that the condition " for each non zero left ideal $I$ of $R$, $I$
contains a uniform left ideal of $R$ " in the Corollary 2.13 is necessary.
In order to see this, consider de following example

\bigskip

\textbf{Example 2.14}. Let $R$ which is an Ore domain to the right but not
to the left. See [15 p 53] for examples of rings of this sort. It is proved
in [7, \ p 486] that there are no $\chi \left( R\right) $-cocritical left $R$%
-modules. As $R$ is a domain, then $R$ is non singular. Thus $\chi \left(
R\right) =\tau _{g}$. Moreover $R$ is a prime ring and $R$ clearly satisfies
ACC on left annihilators. On the other hand we know that if $U$ is a uniform
$\tau _{g}$-torsion free module, then $U$ is $\tau _{g}$-cocritical. But
this is not possible. Thus there are no $\tau _{g}$-torsion free uniform
modules in $R$-mod. Hence $\mathcal{E}_{\tau _{g}}\left( R\right) =\emptyset
$. As $R$ is a domain, then $Spec^{Min}\left( R\right) =\left\{ 0\right\} $.
Thus there are no bijective correspondence between $\mathcal{E}_{\tau
_{g}}\left( R\right) \,$and $Spec^{Min}\left( R\right) $.

\bigskip\

\textbf{Corollary 2.15}. Let $R$ be a semiprime left Goldie ring, then there
is a bijective correspondence between a complete set of representatives of
isomorphism classes of indecomposable non singular injective $R$-modules and
the set of minimal prime ideals of $R$.

\bigskip

\textbf{Proof}. As $R$ is a left Goldie ring, then $R$ is a Goldie module.
So by the Corollary 2.13 we have the result.

\bigskip

Notice that if $R$ is as the Corollary 2.15, then by Theorem 2.7 and
Corollary 2.13 $Spec_{\tau _{g}}\left( R\right) =\left\{ I\mid I\text{ is a }%
\tau _{g}\text{-pure prime ideal of }R\right\} =Spec^{Min}\left( R\right) $.
Thus $R$ has local Gabriel correspondence with respec to $\chi \left(
R\right) =\tau _{g}$

\bigskip

\textbf{Lemma 2.16}. Let $M\in R$-Mod be projective in $\sigma \left[ M%
\right] $ and a semiprime module. Suppose that $M$ satisfies ACC on
annihilators. If $P$ is a minimal prime in $M$ and $Ann_{M}^{r}\left(
P\right) =L$, then $P=Ann_{M}\left( L^{\prime }\right) $ for all $0\neq
L^{\prime }\subseteq L$.

\bigskip

\textbf{Proof.} By Theorem 2.2, $P$ is annihilator submodule, then by
Proposition 1.18 we have that $P=Ann_{M}^{r}\left( Ann_{M}^{r}\left(
P\right) \right) =Ann_{M}\left( L\right) $. Now let $0\neq L^{\prime
}\subseteq L$ and $Ann_{M}\left( L^{\prime }\right) =K$. Thus $%
P=Ann_{M}\left( L\right) \subseteq Ann_{M}\left( L^{\prime }\right) =K$.

Suppose $P\varsubsetneq K$. As $P$ is prime in $M$ and $Ann_{M}\left(
K\right) _{M}K=0\subseteq P$, then $Ann_{M}\left( K\right) \subseteq P$. So $%
Ann_{M}\left( K\right) \subseteq K$. Thus $Ann_{M}\left( K\right)
_{K}Ann_{M}\left( K\right) \subseteq Ann_{M}\left( K\right) _{M}K=0$. Since $%
M$ is semiprime, then $Ann_{M}\left( K\right) =0$. By Proposition 1.18 we
have that $Ann_{M}\left( Ann_{M}\left( K\right) \right) =K$. But $%
Ann_{M}\left( Ann_{M}\left( K\right) \right) =Ann_{M}\left( 0\right) =M$.
Thus $M=K=Ann_{M}\left( L^{\prime }\right) $. So by Proposition 1.9 $%
L^{\prime }=0$ a contradiction. Therefore $P=K=Ann_{M}\left( L^{\prime
}\right) $.

Note that if $M$ is as in Lemma 2.16 and $Ann_{M}^{r}\left( P\right) =L$,
then $Ass_{M}\left( L^{\prime }\right) =\left\{ P\right\} $ for all $0\neq
L^{\prime }\subseteq L$. Hence $\left\{ P\right\} =Ass_{M}\left( N\right)
=Ass_{M}\left( \widehat{L}\right) $ for all $0\neq N\subseteq \widehat{L}$.

\bigskip

\textbf{Proposition 2.17}. Let $M\in R$-Mod be projective in $\sigma \left[ M%
\right] $ and a semiprime module. Suppose that $M$ satisfies ACC on
annihilators. If $P_{1},P_{2},...,P_{n}$ are the minimal primes in $M$, then
$\left\{ N_{1},N_{2},...N_{n}\right\} $ is an independent family, where $%
N_{i}=Ann_{M}^{r}\left( P_{i}\right) $ for $1\leq i\leq n$.

$\bigskip $

\textbf{Proof.} By induction. If $n=1$, then we have the result. Suppose
that $\left\{ N_{1},N_{2},...N_{n-1}\right\} $ is an independent family. If $%
\left( N_{1}\oplus N_{2}\oplus ...\oplus N_{n-1}\right) \cap N_{n}\neq 0$,
then there are $x\in N_{n}$ \ and $1\leq i\leq n-1$, such that $Rx\cong
N_{i}^{\prime }\subseteq N_{i}$. By Lemma 2.16 we have that $%
P_{n}=Ann_{M}^{r}\left( Rx\right) =Ann_{M}^{r}\left( N_{i}^{\prime }\right)
=Ann_{M}^{r}\left( N_{i}\right) =P_{i}$ a contradiction. Therefore $\left(
N_{1}\oplus N_{2}\oplus ...\oplus N_{n-1}\right) \cap N_{n}=0$. So \ \ \ \ \
\ \ $\left\{ N_{1},N_{2},...N_{n}\right\} $ is a independent family$.$\ \ \
\ \ \ \

\bigskip\

\textbf{Theorem 2.18}. Let $M\in R$-Mod be projective in $\sigma \left[ M%
\right] $ and a semiprime module. Suppose that $M$ satisfies ACC on
annihilators and for each $0\neq N\subseteq M$, $N$ contains a uniform
submodule. If $P_{1},P_{2},...,P_{n}$ are the minimal primes in $M$, then $%
\widehat{N_{1}}\oplus \widehat{N_{2}}\oplus ...\oplus \widehat{N_{n}}=%
\widehat{M}$ where $N_{i}=Ann_{M}^{r}\left( P_{i}\right) $ for $1\leq i\leq
n $.

\bigskip

\textbf{Proof.} By Proposition 2.17 we have that $\left\{
N_{1},N_{2},...N_{n}\right\} $ is an independent family. We claim that $%
N_{1}\oplus N_{2}\oplus ...\oplus N_{n}\subseteq _{ess}M$. In fact let $%
0\neq L\subseteq M$. By hypothesis there exists $U$ a uniform module shut
that $U\subseteq L$. By Proposition 1.11 and Theorem 2.2 we have that $%
Ann_{M}\left( U\right) $ is minimal prime in $M$. So $Ann_{M}\left( U\right)
=P_{j}$ for some $1\leq j\leq n$. Hence $\left( P_{j}\right) _{M}U=0$. Thus $%
U\subseteq Ann_{M}^{r}\left( P_{j}\right) =N_{j}$. Whence $\left(
N_{1}\oplus N_{2}\oplus ...\oplus N_{n}\right) \cap L\neq 0$. Therefore $%
\widehat{N_{1}}\oplus \widehat{N_{2}}\oplus ...\oplus \widehat{N_{n}}=%
\widehat{M}$.

\bigskip

\textbf{Corollary 2.19. }Let $M\in R$-Mod be projective in $\sigma \left[ M%
\right] $ and a semiprime module. Suppose that $M$ is Goldie Module and $%
P_{1},P_{2},...,P_{n}$ are the minimal primes in $M$, then $\widehat{N_{1}}%
\oplus \widehat{N_{2}}\oplus ...\oplus \widehat{N_{n}}=\widehat{M}$ where $%
N_{i}=Ann_{M}^{r}\left( P_{i}\right) $ for $1\leq i\leq n$.

\bigskip

\textbf{Proof}. As $M$ is Goldie module, then $M$ has finite uniform
dimension and $M$ satisfies ACC on annihilators. So by the Theorem 2.18 we
have the result.

\bigskip

\bigskip

\textbf{Theorem 2.20.} Let $M\in R$-Mod be projective in $\sigma \left[ M%
\right] $ and a semiprime module. Suppose that $M$ is Goldie Module and $%
P_{1},P_{2},...,P_{n}$ are the minimal primes in $M$ submodules. If $%
N_{i}=Ann_{M}^{r}\left( P_{i}\right) $ for $1\leq i\leq n$, then there
exists $E_{1}$, $E_{2}$, ... , $E_{n}$ uniform injective modules such that $%
\widehat{M}\cong E_{1}^{k_{1}}\oplus E_{2}^{k_{2}}\oplus ...\oplus
E_{n}^{k_{n}}$ and $Ass_{M}\left( E_{i}\right) =\left\{ P_{i}\right\} $.

\bigskip

\bigskip

\textbf{Proof. }Since $N_{i}$ has finite uniform dimension, then there are $%
U_{i_{1}}$, $U_{i_{2}}$, ... , $U_{i_{k_{i}}}$ uniform submodules of $N_{i}$%
, such that $U_{i_{1}}\oplus U_{i_{2}}\oplus ...\oplus
U_{i_{k_{i}}}\subseteq _{ess}N_{i}$. Thus $\widehat{U_{i_{1}}}\oplus
\widehat{U_{i_{2}}}\oplus ...\oplus \widehat{U_{i_{k_{i}}}}=\widehat{N_{i}}$%
. Now by Lemma 2.16, we have that $Ann_{M}\left( U_{i_{j}}\right) =P_{i}$
for all $1\leq j\leq k_{i}$. Thus $Ass_{M}\left( U_{i_{j}}\right)
=Ass_{M}\left( \widehat{U_{i_{j}}}\right) =\left\{ P_{i}\right\} $ for all $%
1\leq j\leq k_{i}$. Hence by Theorem\textbf{\ }2.7, we have that $\widehat{%
U_{i_{1}}}\cong \widehat{U_{i_{2}}}\cong ....\cong \widehat{U_{i_{k_{i}}}}$%
.\ So we can denote $E_{i}=\widehat{U_{i_{1}}}$. Hence $\widehat{N_{i}}\cong
E_{i}^{k_{i}}$ for $1\leq i\leq n$. Therefore by Theorem 2.18 we have that $%
\widehat{M}\cong E_{1}^{k_{1}}\oplus E_{2}^{k_{2}}\oplus ...\oplus
E_{n}^{k_{n}}$.

\bigskip

Notice that by Theorem 2.2 $E_{i}$ $\ncong $ $E_{j}$ for $i\neq j$.

\bigskip

\textbf{Proposition 2.21}. Let $M\in R$-Mod be projective in $\sigma \left[ M%
\right] $ and a semiprime module. Suppose that $M$ is a Goldie Module and $%
P_{1},P_{2},...,P_{n}$ are the minimal primes in $M$, then $\widehat{N_{i}}$
is a fully invariant submodule of $\ \widehat{M}$ where $N_{i}=Ann_{M}^{r}%
\left( P_{i}\right) $ for $1\leq i\leq n$.

\bigskip

\textbf{Proof}. We claim that if $i\neq j$, then $Hom_{R}\left( \widehat{%
N_{i}},\text{ }\widehat{N_{j}}\right) =0$. In fact by Theorem 2.20 we have
that $\widehat{U_{i_{1}}}\oplus \widehat{U_{i_{2}}}\oplus ...\oplus \widehat{%
U_{i_{k_{i}}}}=\widehat{N_{i}}$ where $U_{i_{1}}$, $U_{i_{2}}$, ... , $%
U_{i_{k_{i}}}$ are uniform submodules of $N_{i}$ \ analogously $\ \ \widehat{%
U_{j_{1}}}\oplus \widehat{U_{j_{2}}}\oplus ...\oplus \widehat{U_{j_{k_{j}}}}=%
\widehat{N_{j}}$ where $U_{j_{1}}$, $U_{i_{j2}}$, ... , $U_{j_{k_{j}}}$ are
uniform submodules of $N_{j}$. Let $0\neq f\in Hom_{R}\left( \widehat{N_{i}},%
\text{ }\widehat{N_{j}}\right) $, then there exist $i_{r}$ and $j_{t}$ such
that the restriction morphism. $f_{\mid \widehat{U_{i_{r}}}}:$ $\widehat{%
U_{i_{r}}}\rightarrow $ $\widehat{U_{j_{t}}}$ is non zero. By [4,
Proposition 3.4] $M$ is non $M$-singular. Hence $\widehat{M}$ \ is non $M$%
-singular. Since $\widehat{U_{i_{r}}}$ is uniform submodule of $\widehat{M}$%
, then $f_{\mid \widehat{U_{i_{r}}}}$ is a monomorphism. So $Ass_{M}\left(
\widehat{U_{i_{r}}}\right) =Ass_{M}\left( \widehat{U_{j_{t}}}\right) $. But $%
Ass_{M}\left( \widehat{U_{i_{r}}}\right) =Ass_{M}\left( \widehat{N_{i}}%
\right) =P_{i}$ and $Ass_{M}\left( \widehat{U_{j_{t}}}\right) =Ass_{M}\left(
\widehat{N_{j}}\right) =P_{j}$ a contradiction. Therefore $Hom_{R}\left(
\widehat{N_{i}},\text{ }\widehat{N_{j}}\right) =0$. Now let $g:$ $\widehat{M}%
\rightarrow \widehat{M}$ be a morphism. By Theorem 2.18 we have $g\left(
\widehat{N_{i}}\right) $ $\subseteq $ $\widehat{N_{1}}\oplus \widehat{N_{2}}%
\oplus ...\oplus \widehat{N_{n}}$. Since $Hom_{R}\left( \widehat{N_{i}},%
\text{ }\widehat{N_{j}}\right) =0$ for $i\neq j$, then $g\left( \widehat{%
N_{i}}\right) \subseteq \widehat{N_{i}}$. Thus $\widehat{N_{i}}$ is a fully
invariant submodule of $\ \widehat{M}$.

\bigskip

\bigskip

\bigskip \textbf{Proposition 2.22}. Let $M\in R$-Mod be projective in $%
\sigma \left[ M\right] $ and a semiprime module. Suppose that $M$ satisfies
ACC on annihilators and for each $0\neq N\subseteq M$, $N$ contains a
uniform submodule. If $P_{1},P_{2},...,P_{n}$ are the minimal primes in $M$
and $N_{i}=Ann_{M}^{r}\left( P_{i}\right) $ has finite uniform dimension for
each $1\leq i\leq n$, then $M$ has finite uniform dimension.

\bigskip

\textbf{Proof.} By Theorem 2.18 we obtain that $N_{1}\oplus N_{2}\oplus
...\oplus N_{n}\subseteq _{ess}M$. Now since $N_{i}$ has finite uniform
dimension for each $1\leq i\leq n$, then $M$ has finite uniform dimension.

\bigskip

\textbf{Theorem 2.23.} Let $M\in R$-Mod be projective in $\sigma \left[ M%
\right] $ and a semiprime module. Suppose that $M$ satisfies ACC on
annihilators and for each $0\neq N\subseteq M$, $N$ contains a uniform
submodule. If $P_{1},P_{2},...,P_{n}$ are the minimal primes in $M$ then the
following conditions are equivalent.

$i)$ $M$ is a Goldie Module.

$ii)$ $N_{i}=Ann_{M}^{r}\left( P_{i}\right) $ has finite uniform dimension
for each $1\leq i\leq n$.

\bigskip

Proof $i)\Rightarrow ii)$ As $M$ is left Goldie module, then $M$ has finite
uniform dimension. Hence we have the result.

$ii)\Rightarrow i)$ \ By Proposition 2.22 we have that $M$ has finite
uniform dimension. Thus $M$ is left Goldie module.

\bigskip

\textbf{Corollary 2.24.} Let $R$ be a semiprime ring such that $R$ satisfies
ACC on left annihilators. Suppose that for each non zero left ideal $%
I\subseteq R$, $I$ contains a uniform left ideal. If $P_{1},P_{2},...,P_{n}$
are the minimal prime ideals of $R$, then the following conditions are
equivalent.

\bigskip

$i)$ $R$ is a left Goldie ring.

$ii)$ $Ann_{R}\left( P_{i}\right) =N_{i}$ has finite uniform dimension for
each $1\leq i\leq n$.

\bigskip

Proof. It follows from Theorem 2.23

\bigskip

\textbf{Lemma 2.25.} Let $M\in R$-Mod be projective in $\sigma \left[ M%
\right] $ and a semiprime module. Suppose that $M$ satisfies ACC on
annihilators. If $P_{1},P_{2},...,P_{n}$ are the minimal primes in $M$, then
$Ann_{M}^{r}\left( P_{i}\right) =\underset{i\neq j}{\cap }P_{j}$ for every $%
1\leq i\leq n$.

\bigskip

\textbf{Proof}. As each $P_{i}$ is a fully invariant submodule of $M$, then $%
\left( P_{i}\right) _{M}\left( \underset{i\neq j}{\cap }P_{j}\right)
\subseteq P_{i}\cap \left( \underset{i\neq j}{\cap }P_{j}\right) =\cap
_{j=1}^{n}P_{j}$. By Theorem 2.2 we have that $\cap _{j=1}^{n}P_{j}=0$. Thus
$\left( P_{i}\right) _{M}\left( \underset{i\neq j}{\cap }P_{j}\right) =0$.
So $\underset{i\neq j}{\cap }P_{j}\subseteq Ann_{M}^{r}\left( P_{i}\right) $%
. Now let $K\subseteq M$ such that $\left( P_{i}\right) _{M}K=0$, then $%
\left( P_{i}\right) _{M}K\subseteq P_{j}$ \ for all $i\neq j$. As $P_{j}$ is
prime in $M$, then $P_{i}\subseteq P_{j}$ or $K\subseteq P_{j}$. But $%
P_{i}\subseteq P_{j}$ is not possible. Thus $K\subseteq P_{j}$ for all $%
j\neq i$. Hence $K\subseteq \underset{i\neq j}{\cap }P_{j}$. Thus $%
Ann_{M}^{r}\left( P_{i}\right) =\underset{i\neq j}{\cap }P_{j}$

\bigskip

Notice that by Proposition 1.16 we have that $Ann_{M}\left( P_{i}\right)
=Ann_{M}^{r}\left( P_{i}\right) $. Thus $Ann_{M}\left( P_{i}\right) =%
\underset{i\neq j}{\cap }P_{j}$ for every $1\leq i\leq n$.

\bigskip

\textbf{Corollary 2.26}. Let $M\in R$-Mod be projective in $\sigma \left[ M%
\right] $ and a semiprime module. Suppose that $M$ satisfies ACC on
annihilators and for each $0\neq N\subseteq M$, $N$ contains a uniform
submodule. If $P_{1},P_{2},...,P_{n}$ are the minimal primes in $M$ and $%
P_{i}$ has finite uniform dimension for each $1\leq i\leq n$, then $M$ has
finite uniform dimension.

\bigskip

\textbf{Proof}. By Lemma 2.25 $Ann_{M}^{r}\left( P_{i}\right) =\underset{%
i\neq j}{\cap }P_{j}$. Hence $Ann_{M}^{r}\left( P_{i}\right) $ has finite
uniform dimension. So by Proposition 2.23 we have that $M$ has finite
uniform dimension.

\bigskip

\textbf{Theorem 2.27}. Let $R$ be a semiprime ring such that $R$ satisfies
ACC on left annihilators. Suppose that for each non zero left ideal $%
I\subseteq R$, $I$ contains a uniform left ideal. If $P_{1},P_{2},...,P_{n}$
are the minimal prime ideals of $R$, then the following conditions are
equivalent.

\bigskip

$i)$ $R$ is left Goldie ring.

$ii)$ $P_{i}$ has finite uniform dimension for each $1\leq i\leq n$.

\bigskip

\textbf{Proof}. $i)\Rightarrow ii)$ It is clear.

$ii)\Rightarrow i)$ By Corollary 2.26 we have that $R$ has finite uniform
dimension. So $R$ is left Goldie ring.

\bigskip

Finally we obtain the following result which extends the result given in [
Bo Stemtrom Lemma 2.5]

\bigskip

\textbf{Definition 2.28}. Let $M\in R$-Mod. We say that a submodule $N$ of $%
M $ is a nilpotent submodule if $N^{k}=0$ for some positive integer $k$.

\bigskip

\bigskip \textbf{Proposition 2.29. }Let $M$ be projective in $\sigma \left[ M%
\right] $ . If $M$ satisfies ACC on annihilators, then $\mathcal{Z}\left(
M\right) $ is a nilpotent submodule.

$\bigskip $

\textbf{Proof}. We consider the descending chain $\mathcal{Z}\left( M\right)
\supseteq \mathcal{Z}\left( M\right) ^{2}\supseteq ....$We suppose that
\bigskip $\mathcal{Z}\left( M\right) ^{n}\neq 0$ for all $n\geq 1$.So we
have the ascending chain $Ann_{M}\left( \mathcal{Z}\left( M\right) \right)
\subseteq Ann_{M}\left( \mathcal{Z}\left( M\right) ^{2}\right) \subseteq ...$%
Since $M$ satisfies ACC on annihilators, there exists $n>0$ such that $%
Ann_{M}\left( \mathcal{Z}\left( M\right) ^{n}\right) =Ann_{M}\left( \mathcal{%
Z}\left( M\right) ^{n+1}\right) =....$

As $\mathcal{Z}\left( M\right) ^{n+2}\neq 0$, then by [1, Proposition 5.6]
we have that $\left[ \mathcal{Z}\left( M\right) _{M}\mathcal{Z}\left(
M\right) \right] _{M}\mathcal{Z}\left( M\right) ^{n}\neq 0$.

Since $\mathcal{Z}\left( M\right) _{M}\mathcal{Z}\left( M\right) =\sum
\left\{ f(\mathcal{Z}\left( M\right) )\mid f:M\rightarrow \mathcal{Z}\left(
M\right) \text{ }\right\} $, then\

$\left[ \sum \left\{ f(\mathcal{Z}\left( M\right) )\mid f:M\rightarrow
\mathcal{Z}\left( M\right) \text{ }\right\} \right] _{M}\mathcal{Z}\left(
M\right) ^{n}\neq 0$. So by [3, Proposition1.3 ] there exits $f:M\rightarrow
\mathcal{Z}\left( M\right) $ such that $f\left( \mathcal{Z}\left( M\right)
\right) _{M}\mathcal{Z}\left( M\right) ^{n}\neq 0$. Thus $f\left( M\right)
_{M}\mathcal{Z}\left( M\right) ^{n}\neq 0$

Now consider the set

\begin{equation*}
\Gamma =\left\{ \ker f\mid f:M\rightarrow \mathcal{Z}\left( M\right) \text{
and }f\left( M\right) _{M}\mathcal{Z}\left( M\right) ^{n}\neq 0\right\}
\end{equation*}

By hypothesis $\Gamma $ has maximal elements. Let $f:M\rightarrow \mathcal{Z}%
\left( M\right) $ such that $\ker f$ is a maximal element in $\Gamma $.

\bigskip

\bigskip

We claim that $h\left( f(M)\right) _{M}\mathcal{Z}\left( M\right) ^{n}=0$
for all $h:M\rightarrow \mathcal{Z}\left( M\right) $. In fact let $%
h:M\rightarrow \mathcal{Z}\left( M\right) $ be a morphism. By [5, Lemma 2.7]
$\ker h\subseteq _{ess}M$. So $\ker h\cap f\left( M\right) \neq 0$. Thus
there exists $0\neq f\left( m\right) $ such that $h\left( f\left( m\right)
\right) =0$. So $\ker f\varsubsetneq \ker h\circ f$. If $h\left( f(M)\right)
_{M}\mathcal{Z}\left( M\right) ^{n}\neq 0$, then $h\circ f\in \Gamma $. But $%
\ker f$ is a maximal element in $\Gamma $ and $\ker f\varsubsetneq \ker
h\circ f$, then $h\circ f\notin \Gamma $ a contradiction. Therefore $h\left(
f(M)\right) _{M}\mathcal{Z}\left( M\right) ^{n}=0$.

On the other hand by [1, Proposition 5.6] we have that .

$f\left( M\right) _{M}\mathcal{Z}\left( M\right) ^{n+1}=\left[ f\left(
M\right) _{M}\mathcal{Z}\left( M\right) \right] _{M}\mathcal{Z}\left(
M\right) ^{n}=$

$\left[ \sum \left\{ h\left( f\left( M\right) \right) \mid h:M\rightarrow
\mathcal{Z}\left( M\right) \right\} \right] _{M}\mathcal{Z}\left( M\right)
^{n}$. By [3, Proposition 1.3] we have that $\left[ \sum \left\{ h\left(
f\left( M\right) \right) \mid h:M\rightarrow \mathcal{Z}\left( M\right)
\right\} \right] _{M}\mathcal{Z}\left( M\right) ^{n}=$

$\sum \left\{ h\left( f\left( M\right) \right) _{M}\mathcal{Z}\left(
M\right) ^{n}\mid h:M\rightarrow \mathcal{Z}\left( M\right) \right\} =0$.

Thus $f\left( M\right) _{M}\mathcal{Z}\left( M\right) ^{n+1}=0$.

Hence $f\left( M\right) \subseteq Ann_{M}\left( \mathcal{Z}\left( M\right)
^{n+1}\right) =Ann_{M}\left( \mathcal{Z}\left( M\right) ^{n}\right) $. So $%
f\left( M\right) _{M}\mathcal{Z}\left( M\right) ^{n}=0$ a contradiction.

\bigskip

\textbf{Corollary 2.30. }Let $M\in R$- Mod be projective in $\sigma \left[ M%
\right] $ and $S=End_{R}\left( M\right) $. If $M$ is retractable and
satisfies ACC on annihilators, then $\mathcal{Z}_{r}\left( S\right) $ is
nilpotent. Where $\mathcal{Z}_{r}\left( S\right) $ is the right singular
ideal of $S$.

\bigskip

\textbf{Proof}. Notice that if, $N$ and $L$ are submodules of $M$ then,

$Hom_{R}\left( M,L\right) Hom_{R}\left( M,N\right) \subseteq Hom_{R}\left(
M,N_{M}L\right) $. Now we consider the ideal $\Delta =\left\{ f\in S\mid
\ker f\subseteq _{ess}M\right\} $. We claim that $\mathcal{Z}_{r}\left(
S\right) \subseteq \Delta $, in fact let $\alpha \in \mathcal{Z}_{r}\left(
S\right) $, then there exists $I$ essential right deal of $S$ such that $%
\alpha I=0$. Thus $\alpha \circ g=0$ for all $\ g\in I$. Hence $\alpha
\left( \sum\limits_{g\in I}g\left( M\right) \right) =0$. Now let $0\neq
N\subseteq M$ . As $M$ is retractable, then there exits $\rho :M\rightarrow
N $ a non zero morphism. Thus $\rho S\cap I\neq 0$. So there exists $h\in S$
such that $0\neq \rho \circ h\in I$. Hence $0\neq \left( \rho \circ h\right)
\left( M\right) \subseteq N\cap \sum\limits_{g\in I}g\left( M\right) $. Thus
$\sum\limits_{g\in I}g\left( M\right) \subseteq _{ess}M$. So $\ker \alpha
\subseteq _{ess}M$. Hence $\mathcal{Z}_{r}\left( S\right) \subseteq \Delta $%
. So $\mathcal{Z}_{r}\left( S\right) M\subseteq \Delta M\subseteq \mathcal{Z}%
\left( M\right) $. Moreover $Hom_{R}\left( M,\Delta M\right) \subseteq
\left( M,\mathcal{Z}\left( M\right) \right) $. On the other hand we have
that $\Delta \subseteq Hom_{R}\left( M,\Delta M\right) $. By Proposition
2.29 there exists $n>0$ such that $\mathcal{Z}\left( M\right) ^{n}=0$. Thus $%
\Delta ^{n}\subseteq Hom_{R}\left( M,\Delta M\right) ^{n}\subseteq \left( M,%
\mathcal{Z}\left( M\right) \right) ^{n}\subseteq \left( M,\mathcal{Z}\left(
M\right) ^{n}\right) =0$. So $\Delta $ is nilpotent. Since $\mathcal{Z}%
_{r}\left( S\right) \subseteq \Delta $, then $\mathcal{Z}_{r}\left( S\right)
$ is nilpotent.

\bigskip

In the Corollary 2.30 we have that $\Delta $ is nilpotent. Notice that to
obtain this resul it is not necessary the condition $M$ is retractable.

\bigskip

\bigskip

\textbf{Corollary 2.31. }Let $M\in R$- Mod be projective in $\sigma \left[ M%
\right] $ and $S=End_{R}\left( M\right) $. If $I$ is an ideal of $S$ such
that $\cap _{f\in I}\ker f\subseteq _{ess}M$, then $I$ is nilpotent.

\bigskip

\textbf{Proof.} As $\cap _{f\in I}\ker f\subseteq _{ess}M$, then $\ker
f\subseteq _{ess}M$ for all $\ f\in I$. Thus $I\subseteq \Delta $. By
Corollary 2.30 we have that $\Delta $ is nilpotent. So $I$ is nilpotent

\bigskip

\bigskip

\textbf{Definition 2.32.} Let $M\in R$. A submodule $N$ of $M$ is $T_{M}$%
-nilpotent in case every sequence $\left\{ f_{1},f_{2},...,f_{n},...\right\}
$ in $Hom_{R}\left( M,N\right) $ and any $a\in N$, there exists $n\geq 1$
such that $f_{n}f_{n-1}...f_{1}\left( a\right) =0$.

\bigskip

Notice that if $I$ is a left ideal of a ring $R$, then $I$ is left $T$%
-nilpotent in the usual sense if and only if $I$ \ is $T_{R}$-nilpotente.

\bigskip

\textbf{Proposition 2.33}. Let $M\in R$- Mod be projective in $\sigma \left[
M\right] $ and retractable. Suppose that $M$ satisfies ACC on annihilators.
If \ $N\subseteq M$ is $T_{M}$-nilpotent, then $N$ is nilpotent.

\bigskip

\textbf{Proof}. Let $N\subseteq M$ be $T_{M}$-nilpotent. Consider the chain $%
N\supseteq N^{2}\supseteq N^{3}\supseteq ...$. Then we have the chain $%
Ann_{M}\left( N\right) \subseteq Ann_{M}\left( N^{2}\right) \subseteq
Ann_{M}\left( N^{3}\right) \subseteq ...$. Since $M$ satisfies ACC on
annihilators, then there exists $k\geq 1$ such that $Ann_{M}\left(
N^{K}\right) =Ann_{M}\left( N^{k+1}\right) =Ann_{M}\left( N^{k+2}\right)
.....$ Let $L=N^{k}$. If $L^{2}=L_{M}L\neq 0$, then there exists $%
f:M\rightarrow L$ and $a\in L\subseteq N$ such that $f\left( a\right) \neq 0$%
. Hence $\left( Ra\right) _{M}L\neq 0$. If $\left( Ra\right) _{M}\left(
L_{M}L\right) =0$, then $Ra\subseteq Ann_{M}\left( L^{2}\right)
=Ann_{M}\left( L\right) $. So $Ra_{M}L=0$ a contradiction. Thus $\left(
Ra\right) _{M}\left( L_{M}L\right) \neq 0$. By [1, M-injective 5.6] $0\neq
\left( Ra\right) _{M}\left( L_{M}L\right) =\left( \left( Ra\right)
_{M}L\right) _{M}L$. As $\left( Ra\right) _{M}L=\sum \left\{ f\left(
Ra\right) \mid f:M\rightarrow L\right\} $, then by [3, Proposition 1.3]
there exists $f_{1}:M\rightarrow L$ such that $f_{1}\left( Ra\right)
_{M}L\neq 0$. So $\left( Rf_{1}\left( a\right) \right) _{M}L\neq 0$. Newly $%
\left( Rf_{1}\left( a\right) \right) _{M}\left( L_{M}L\right) \neq 0$. Hence
$\left( Rf_{1}\left( a\right) _{M}L\right) _{M}L\neq 0$. So there exists $%
f_{2}:M\rightarrow L$ such that $f_{2}\left( Rf_{1}\left( a\right) \right)
_{M}L\neq 0$. Hence $Rf_{2}\left( f_{1}\left( a\right) \right) _{M}L\neq 0$.
Continuing in this way, we have that $Rf_{n}.f_{n-1}...f_{1}\left( a\right)
_{M}L\neq 0$ for all $n\geq 1$. Hence $f_{n}.f_{n-1}...f_{1}\left( a\right)
\neq 0$ for all $n\geq 1$ a contradiction because $N$ is $T_{M}$-nilpotent.
Hence $L^{2}=0$. So $N$ is nilpotent.

\bigskip

\bigskip

\section{\protect\bigskip Continuous Modules with ACC on annihilators}

\bigskip

\bigskip The following definitions were given in [14]. We include here these
definitions for convenience of the reader.

\bigskip

\textbf{Definition 3.1.} Let $M\in R$-Mod. $M$ is called continuous module
if it satisfies the followings conditions:

$C_{1})$ \ Every submodule of $M$ is essential in a direct summand of $M$

$C_{2})$ \ Every submodule of $M$ that is isomorphic to a direct summand of $%
M$ is itself a direct summand.

\bigskip

\textbf{Definition 3.2}. An $R$-module $M$ is $\mathcal{K}$-nonsingular if
for every $f\in End_{R}\left( M\right) $ such that $\ker f\subseteq _{ess}M$
implies $f=0$.

\bigskip

Notice that by [14, Proposition 2.3] we have that if $M$ is non $M$-singular
module, then $M$ is $\mathcal{K}$-nonsingular

\bigskip

\textbf{Definition 3.3}. Let $M\in R$-Mod. The $\mathcal{K}$-singular
submodule of $M$ is defined as $Z^{\mathcal{K}}\left( M\right) =\sum \left\{
f\left( M\right) \mid f\in End_{R}\left( M\right) \text{ and }\ker
f\subseteq _{ess}M\text{ }\right\} $

\bigskip

\textbf{Remark 3.4}. If $M$ is a continuous module and \ $S=End_{R}\left(
M\right) $, then by [10, Proposition 1.25 ] we have that $J\left( S\right)
=\left\{ f\in S\mid \ker f\subseteq _{ess}M\right\} $. Hence $Z^{\mathcal{K}%
}\left( M\right) =J\left( S\right) M$. On the other hand it is clear that $%
Z^{\mathcal{K}}\left( M\right) \subseteq \mathcal{Z}\left( M\right) $.
Therefore $J\left( S\right) M\subseteq \mathcal{Z}\left( M\right) $.

\bigskip

\textbf{Proposition 3.5.} Let $M\in R$- Mod be projective in $\sigma \left[ M%
\right] $ and $S=End_{R}\left( M\right) $. Suppose that $M$ is a continuous
module . If $M$ satisfies ACC on annihilators, then $J(S)$ is nilpotent

\bigskip

\textbf{Proof.} As $M$ is a continuous module, then by Remark 3.4 we have
that $J\left( S\right) =\left\{ f\in S\mid \ker f\subseteq _{ess}M\right\} $%
. By Corollary 2.30 $J\left( S\right) =\Delta $ is nilpotent.

\bigskip

\textbf{Theorem 3.6.} Let $M\in R$- Mod be projective in $\sigma \left[ M%
\right] $ and $S=End_{R}\left( M\right) $. Suppose that $M$ is a continuous,
retractable, non $M$-singular module and satisfies ACC on annihilators. Then
$M$ is a semiprime Goldie module.

\bigskip

\textbf{Proof. }Since $M$ is non $M$-singular and continuous, then by Remark
3.4. $J\left( S\right) M=0$. Hence $J\left( S\right) =0$. \ Thus $S$ is a
semiprime ring. We claim that $M$ is semiprime module. In fact let $L$ be a
fully invariant submodule of $M$ such that $L^{2}=L_{M}L=0$. As $%
Hom_{R}\left( M,L\right) Hom_{R}\left( M,L\right) \subseteq Hom_{R}\left(
M,L_{m}L\right) $, then $Hom_{R}\left( M,L\right) Hom_{R}\left( M,L\right)
=0 $. Since $S$ is semiprime ring the $Hom_{R}\left( M,L\right) =0$. Now as $%
M$ is retractable, then $L=0$. Since $M$ is non $M$-singular, then $M$ is $%
\mathcal{K}$-nonsingular, then by [14, Proposition 3.1, Corollary 3.3 and
Theorem 1.5] we have that $\ker f$ is a direct summand of $M$ for all $f\in
S $. As satisfies ACC on annihilators, then $M$ satisfies ACC on direct
summands. Since $M$ \ satisfies the condition $C_{1}$, then every closed
submodule of $M$ is a direct summand . Thus $M$ satisfies ACC on closed
submodules. By [9, Proposition 6.30] $M$ has finite uniform dimension. Hence
$M$ is a semiprime Goldie module.

\bigskip

\textbf{Corollary 3.7}. Let $R$ be a continuous and non singular ring.
Suppose that $R$ satisfies ACC on left annihilators, then $R$ is a semiprime
left Goldie ring.

\bigskip

Notice that in the Corollary 3.7 the inverse is not true in general.
Consider the following example

\bigskip

\textbf{Example 3.8. }Let $R=%
\mathbb{Z}
$ be the ring of integers. We know $%
\mathbb{Z}
$ is a prime Goldie ring, but $%
\mathbb{Z}
$ is not a continuous ring.

\noindent


\bigskip

\begin{thebibliography}{99}
\bibitem{} Beachy, J. (2002). M-injective modules and prime M-ideals,
\textit{Comm. Algebra} 30(10):4639-4676.

\bibitem{} Bican, L., Jambor, P., Kepka, T., Nemec, P. (1980). Prime and
coprime modules, \textit{Fundamenta Matematicae} 107:33-44

\bibitem{} Castro, J., R\'{\i}os, J. (2012). Prime submodules and local
Gabriel correspondence in $\sigma \left[ M\right] $, \textit{Comm. Algebra}.
40(1):213-232.

\bibitem{} Castro, J., R\'{\i}os, J. (2014). Krull dimension and classical
Krull dimension of modules. \textit{Comm. Algebra. }42(7):3183-3204

\bibitem{} Castro, J., R\'{\i}os, J. (2015). $M$ tame modules and Gabriel
correpondence, \textit{to appear in Comm. Algebra. }

\bibitem{} Chatters, A. , Hajarnavis, C. (1980) Rings with chain conditions,
Boston, London, Melbourne: Pitman Advanced Publishing Program

\bibitem{} Golan, J. (1986). Torsion Theories, Longman Scientific \&
Technical, Harlow.

\bibitem{} Goodearl, K, Warfield, R ( 2004). An introduction to
Noncommutative Noetherian Rings, Cambridge University Press.

\bibitem{} Lam, T. (1998), Lectures on Modules and Rings, Graduate Texts in
Mathematics., vol 139, New York: Springer-Verlag.

\bibitem{} Nicholson, W. Yousif, M. (2003). Quasi-Frobenius Rings. Cambridge
University Press.

\bibitem{} Raggi, F., Rios, J., Rinc\'{o}n, H., Fern\'{a}ndez-Alonso, R.,
Signoret, C. (2005). Prime and irreducible preradicals, \textit{J. Algebra
Appl}. 4(4):451-466.

\bibitem{} Raggi, F., R\'{\i}os, J., Rinc\'{o}n, H., Fern\'{a}ndez-Alonso,
R. (2009). Basic Preradicals and Main Injective Modules, \textit{J. Algebra
Appl.} 8(1):1-16.

\bibitem{} Raggi, F., R\'{\i}os, J., Rinc\'{o}n, H., Fern\'{a}ndez-Alonso,
R. (2009). Semiprime Preradicals, \textit{Comm. Algebra}. 37(7):2811-2822.

\bibitem{} Rizvi, S., Roman, C. (2007). On $K$-nonsigular Modules and
Applications. \textit{Comm. Algebra}. 35:2960-2982..

\bibitem{} Stenstr\"{o}m, B. (1975). Rings of Quotients, Graduate Texts in
Mathematics, New York: Springer-Verlag.

\bibitem{} Wisbauer, R. (1991). Foundations of Module and Ring Theory,
Gordon and Breach: Reading.

\bibitem{} Wisbauer, R.(1996). Modules and Algebras: Bimodule Structure and
Group Actions on Algebras, England: Addison Wesley Longman Limited.
\end{thebibliography}
\end{document}